\documentclass[epsf]{article}
\usepackage{here}
\usepackage{amsmath}
\usepackage{amssymb}
\usepackage{amsfonts}
\usepackage{theorem}
\usepackage{blkarray}
\usepackage{arydshln}
\usepackage{multirow}
\usepackage[dvips]{graphicx}

\newcommand{\bigzerou}{%
\smash{\lower1.7ex\hbox{\bg 0}}}
\newtheorem{theorem}{Theorem}
\newtheorem{prop}{Proposition}
\newtheorem{defi}{Definition}

\newtheorem{conj}{Conjecture}

\newcommand{\ba}{\begin{eqnarray}}
\newcommand{\ea}{\end{eqnarray}}
\newcommand{\ban}{\begin{eqnarray*}}
\newcommand{\ean}{\end{eqnarray*}}
\newcommand{\no}{\nonumber}

\newcommand{\mapright}[1]{%
\smash{\mathop{%
\hbox to 1.0cm{\rightarrowfill}}\limits^{#1}}}
\newcommand{\mapleft}[1]{%
\smash{\mathop{%
\hbox to 1.3cm{\leftarrowfill}}\limits^{#1}}}
\begin{document}
\title{
\begin{flushright}
  \begin{minipage}[b]{5em}
    \normalsize
    ${}$      \\
  \end{minipage}
\end{flushright}
{\bf  Elliptic Virtual Structure Constants and Generalizations of BCOV-Zinger Formula to Projective Fano Hypersurfaces}
\author{Masao Jinzenji${}^{(1)}$, Ken Kuwata${}^{(2)}$ \\
\\
${}^{(1)}$\it Department of Mathematics, Okayama University \\
\it  Okayama, 700-8530, Japan\\
\\
${}^{(2)}$\it Department of General Education\\
\it  National Institute of Technology, Kagawa College \\
\it  Chokushi, Takamatsu, 761-8058, Japan\\
\\
\it e-mail address: \it\hspace{0.3cm}${}^{(1)}$ pcj70e4e@okayama-u.ac.jp \\
\it\hspace{3.6cm}${}^{(2)}$ kuwata-k@t.kagawa-nct.ac.jp}}
\maketitle

\begin{abstract}
In this paper, we propose a method for computing genus $1$ Gromov-Witten invariants of Calabi-Yau and Fano projective hypersurfaces using the B-model. Our formalism is applicable to both Calabi-Yau and Fano cases. In the Calabi-Yau case, significant cancellation of terms within our formalism occurs, resulting in an alternative representation of the BCOV-Zinger formula for projective Calabi-Yau hypersurfaces.
\end{abstract}
\section{Introduction}

\quad The evaluation of higher genus Gromov-Witten invariants has been one of the most exciting and challenging topics in the study of topological strings and mirror symmetry.

In 1993, Bershadsky, Cecotti, Ooguri, and Vafa evaluated the genus 1 Gromov-Witten invariants of the quintic Calabi-Yau 3-fold and several other Calabi-Yau 3-folds by combining the holomorphic anomaly equation with mirror symmetry \cite{BCOV}. Shortly after their paper appeared, they extended the method introduced in \cite{BCOV} to evaluate higher genus ($g \geq 2$) Gromov-Witten invariants of Calabi-Yau 3-folds \cite{BCOV2}. Efforts to generalize these results to higher-dimensional Calabi-Yau manifolds can be found in \cite{sugiyama}. In \cite{FLY}, Fang, Lu, and Yoshikawa derived the BCOV formula from the perspective of analytic torsion.

On the other hand, the study of higher genus Gromov-Witten invariants of Fano manifolds began after the publication of the seminal paper \cite{KM}, which derived the associativity equation of the quantum cohomology ring by analyzing the Chow ring of the moduli space of genus 0 stable curves with four marked points. This breakthrough enabled the evaluation of genus 0 Gromov-Witten invariants of Fano manifolds. In 1997, Getzler generalized the associativity equation to genus 1 Gromov-Witten invariants by analyzing the Chow ring of the moduli space of genus 1 stable curves with four marked points, leading to what is now known as Getzler's equation \cite{G}. In the same year, Eguchi, Hori, and Xiong proposed the Virasoro conjecture, which, in combination with the topological recursion relation, has the potential to evaluate Gromov-Witten invariants at all genera for Fano manifolds \cite{EHX}. These results rely on differential equations governing the generating function of Gromov-Witten invariants.

There is another approach to evaluating Gromov-Witten invariants. In 1994, Kontsevich directly computed genus 0 Gromov-Witten invariants of projective hypersurfaces by combining the moduli space of stable maps with the Bott residue formula \cite{kont}. Although the computational process of this method is highly complex, it offers valuable insights into the geometric study of Gromov-Witten invariants. In 1997, Graber and Pandharipande extended this method \cite{kont} to evaluate higher genus Gromov-Witten invariants for projective spaces \cite{GP}. In their work \cite{GP}, the key point of this generalization is the consideration of the Hodge bundle on the moduli space of higher genus stable curves. As far as we know, the results in \cite{GP} are also applicable to evaluating higher genus local Gromov-Witten invariants of projective spaces \cite{KZ}, but not to compact projective hypersurfaces, including compact Calabi-Yau hypersurfaces. This trend in the field can be summarized by the term "localization."

We summarize the situation of Gromov-Witten invariants research in the 1990s. At the genus 0 level, the mirror symmetry conjecture became more approachable with the development of the localization method. However, this method cannot be applied (at least in the case of compact Calabi-Yau manifolds) to attack the BCOV conjecture in a mathematical sense. On the other hand, the associativity equation is not effective for evaluating genus 0 Gromov-Witten invariants of Calabi-Yau manifolds. Getzler's equation and the Virasoro conjecture are similarly ineffective for evaluating higher genus Gromov-Witten invariants of Calabi-Yau manifolds.

In 2009, Zinger introduced the concept of reduced Gromov-Witten invariants and, by combining it with a type of localization, derived a closed formula for the generating function of genus 1 Gromov-Witten invariants of projective Calabi-Yau hypersurfaces \cite{zinger}. Since then, there has been substantial mathematical literature on higher genus Gromov-Witten invariants of Calabi-Yau 3-folds \cite{EMM, Hu, Ob, Popa}. We refer to the formula proved in \cite{zinger} as the BCOV-Zinger formula, because Zinger's formula for the quintic 3-fold coincides with the one proposed in \cite{BCOV}. Based on extensive experimental computations, we believe his formula is correct, although understanding the logic behind his proof is challenging for us, perhaps due to our background in physics.

{\bf The motivation of this paper is ``to construct our original approach to genus 1 Gromov-Witten invariants of projective hypersurfaces from the perspective of the theory of moduli spaces of quasimaps, which we have developed so far''.}

Let us briefly review our study of classical mirror symmetry for projective hypersurfaces. Let $M_{N}^{k}$ denote a degree $k$ nonsingular hypersurface in ${CP}^{N-1}$. We denote by $h$ the hyperplane class in $H^{1,1}({CP}^{N-1}, \mathbb{C})$. The quantum cohomology ring of $M_{N}^{k}$ is defined by the following multiplication rule:
 \ba
{\cal O}_{h}\cdot{\cal O}_{h^{N-2-m}}=\sum_{d=0}^{\infty}\frac{1}{k}\langle{\cal O}_{h}{\cal O}_{h^{N-2-m}}{\cal O}_{h^{m-1+(N-k)d}}\rangle_{0,d}e^{dt}{\cal O}_{h^{N-1-m-(N-k)d}},
\ea
where $\langle{\cal O}_{h}{\cal O}_{h^{N-2-m}}{\cal O}_{h^{m-1+(N-k)d}}\rangle_{0,d}$ is genus $0$ and degree $d$ Gromov-Witten invariant of $M_{N}^{k}$. 

In \cite{CJ}, we defined $L_{m}^{N,k,d}$ as $\frac{1}{k}\langle{\cal O}_{h}{\cal O}_{h^{N-2-m}}{\cal O}_{h^{m-1+(N-k)d}}\rangle_{0,d}$, which represents the structure constant of the quantum cohomology ring of $M_{N}^{k}$. We also proved recursive formulas describing $L_{m}^{N,k,d}$ in terms of weighted homogeneous polynomials of $L_{m^{\prime}}^{N+1,k,d^{\prime}}$ (with $d^{\prime}\leq d$) up to $d=3$. These recursive formulas take the same form when $N-k\geq 2$. However, for the cases of $N-k=1,0$, these recursive formulas require modification.

We then applied the recursive formulas for the $N-k\geq 2$ cases formally to the $N-k=1,0$ cases and found that the resulting "virtual structure constant" becomes the expansion coefficient of the generalized hypergeometric series used in the mirror computation of Gromov-Witten invariants of $M_{k}^{k}$ and the derivative of the mirror map. We denote the resulting number by $\tilde{L}_{m}^{k,k,d}$ and call it a "virtual structure constant." This marks the beginning of the notion of the "virtual structure constant," which serves as the B-model analog of the structure constant (Gromov-Witten invariant) of the quantum cohomology ring of $M_{k}^{k}$.

In \cite{Jin0, Jin1}, we demonstrated that the virtual structure constant $\tilde{L}_{m}^{k,k,d}$ is given as the intersection number 
$\frac{d}{k}\cdot w({\cal O}_{h^{k-2-m}}{\cal O}_{h^{m-1}})_{0,d}$ on the moduli space of quasimaps from ${CP}^{1}$ with two marked points to ${CP}^{k-1}$, denoted by $\widetilde{Mp}_{0,2}(k,d)$\footnote{For another approach to the moduli space of quasimaps, please refer to \cite{CK}.}. Subsequently, we began referring to $w({\cal O}_{h^{a}}{\cal O}_{h^{b}})_{0,d}$, which is defined as the intersection number on $\widetilde{Mp}_{0,2}(k,d)$ (later generalized to the general hypersurface $M_{N}^{k}$), as the "virtual structure constant." By this term, we mean the B-model analog of Gromov-Witten invariants.

In \cite{JS1}, we constructed the moduli space of quasimaps from ${CP}^{1}$ with $2+n$ marked points to ${CP}^{N-1}$, denoted by $\widetilde{Mp}_{0,2|n}(N,d)$. We also defined the multi-point virtual structure constant $w({\cal O}_{h^{a}}{\cal O}_{h^{b}}|\prod_{p=0}^{N-2}({\cal O}_{h^{p}})^{m_{p}})_{0,d}$ for $M_{N}^{k}$ as the intersection number on 
$\widetilde{Mp}_{0,2|n}(N,d)$. This is expected to serve as the B-model analog of the genus $0$ Gromov-Witten invariant 
$\langle {\cal O}_{h^{a}}{\cal O}_{h^{b}}\prod_{p=0}^{N-2}({\cal O}_{h^{p}})^{m_{p}}\rangle_{0,d}$ for $M_{N}^{k}$, but it does not generically coincide with the corresponding Gromov-Witten invariant.

Next, we defined the "mirror map," which is the generating function of $w({\cal O}_{h^{a}}{\cal O}_{1}|\prod_{p=0}^{N-2}({\cal O}_{h^{p}})^{m_{p}})_{0,d}$ (see the next section for details). We conjectured that the generating function of the multi-point virtual structure constants can be translated into the generating function of Gromov-Witten invariants via a coordinate change of deformation parameters given by the mirror map (see also the next section for details). This conjecture enables us to perform non-trivial mirror computations of Gromov-Witten invariants for Fano projective hypersurfaces, even though the mirror map used in the mirror computation of the structure constants of the small quantum cohomology ring of Fano projective hypersurfaces is trivial. This construction is the starting point of this paper.

At this stage, we discuss the geometric merit of introducing the virtual structure constant or the moduli space of quasimaps. First of all, the geometric structure of the moduli space of quasimaps is far simpler than that of the moduli space of stable maps. This simplicity allows us to present a concise proof of the mirror theorem for the projective hypersurface $M_{N}^{k}$ \cite{Jin2}.

To illustrate the advantage of the moduli space of quasimaps, consider localization computation as an example. In the localization computation of genus $0$ and degree $d$ Gromov-Witten invariants, the result is given as the sum of contributions from tree graphs with $d$ edges. As $d$ increases, the number of such tree graphs grows rapidly. Therefore, the localization computation of Gromov-Witten invariants becomes increasingly difficult as $d$ rises. In contrast, in the localization computation of the virtual structure constant of degree $d$, the result is represented by a multi-variable residue integral corresponding to a contribution from a line graph with $d$ edges. The line graph with $d$ edges is unique, meaning the number of graphs needed to evaluate the genus $0$ virtual structure constant of degree $d$ is always one. Moreover, the multi-variable residue integral representation is more manageable for computer calculations than traditional localization computations, which involve summing the characters of torus actions. As discussed in \cite{JS1}, the mirror transformation from Gromov-Witten invariants to virtual structure constants can be interpreted as the process of "cutting the unnecessary edges of a tree graph into a line graph." In other words, in our context, the mirror transformation from A-model to B-model can be viewed as the "linearization of the moduli space."

In this paper, we present the final results of our attempt to find multi-variable residue integral formulas that correspond to this "linearization of the moduli space" at the level of genus $1$ Gromov-Witten invariants of projective hypersurfaces.

The result obtained by Graber and Pandharipande \cite{GP} represents genus $1$ Gromov-Witten invariants of ${CP}^{n}$ as the sum of contributions from typical tree graphs with one elliptic vertex and one-loop graphs. We first attempted to rewrite the results of genus $1$ Gromov-Witten invariants of ${CP}^{2}$ obtained using the Graber-Pandharipande method into residue integrals. We then translated these into the B-model form using the mirror map introduced in \cite{JS1} (we had already derived the formula for multi-point virtual structure constants). As expected, this process effectively "cut the unnecessary edges of the tree graph with an elliptic vertex and one-loop graph." We observed that the expected residue integral representation could be written as the sum of contributions from graphs of type (i) and type (ii) presented in Section 3 of this paper. However, we also found it necessary to introduce type (iii) and (iv) graphs, also presented in Section 3, to cancel error terms that arise when rewriting localization formulas of the moduli space of stable maps into residue integrals. We informally referred to these error terms as "diagonal anomalies." Controlling this diagonal anomaly was very challenging, but we later discovered a method to manage it in the case of ${CP}^{2}$, which enabled us to derive the final residue integral representation of the elliptic virtual structure constants of ${CP}^{2}$. Ultimately, we succeeded in generalizing the contributions from type (iii) and (iv) graphs to the case of Fano hypersurfaces.

Moreover, we found that these formulas are applicable to Calabi-Yau hypersurfaces $M_{k}^{k}$ as well, yielding a geometrically natural residue integral representation of the BCOV-Zinger formula. In this case, all the contributions from the type (iii) graphs vanish, and our representation clarifies the geometric meaning of the terms appearing in the BCOV formula.

As mentioned earlier, our approach to constructing the residue integral formulas is purely heuristic. Therefore, in this paper, we only present the setup of our mirror computation, the integrands of residue integrals associated with graphs, and the results of numerical tests of our conjectural formulas.

In the original paper \cite{BCOV}, the BCOV formula was derived by considering the "holomorphic anomaly" of the topological conformal field theory. It was obtained by taking the holomorphic limit of the new index of the topological conformal field theory coupled with the complex conjugate of the theory. Thus, a natural question arises: Is our formulation related to the original framework of the BCOV formula? At this point, we believe our results have little connection with the holomorphic anomaly, in the same way that Zinger's derivation has little connection to it. From our perspective, the holomorphic anomaly is just "one way" to derive the B-model generating function of genus $1$ Gromov-Witten invariants of Calabi-Yau manifolds, and the formula we found is simply "another way" to obtain the B-model generating function. However, our results also suggest the possibility of generalizing the BCOV formula to Fano manifolds, and we believe that exploring the connection between our results and the holomorphic anomaly is an interesting avenue of future research.

This paper is organized as follows. In Section 2, we explain the setup of our mirror computation of elliptic Gromov-Witten invariants for projective Fano hypersurfaces, based on the results presented in \cite{JS1}. In Section 3, we present our residue integral representation of elliptic virtual structure constants. In Subsection 3.1, we introduce four types of graphs used in our construction. In Subsection 3.2, we present the integrands of the residue integral representation associated with each type of graph. In Section 4, we test our residue integral representation through numerical computations. In Subsection 4.1, we present results for Fano projective hypersurfaces, observing agreement with known results. In Subsection 4.2, we present results for the Calabi-Yau hypersurface $M_{k}^{k}$, obtaining another representation of the BCOV-Zinger formula. In Appendix A, we compile our numerical results for Fano 3-folds $M_{5}^{k}$ $(k=1,2,3,4)$ in tables.

This paper is written from a physicist's perspective, and the results presented include conjectures that are of interest to mathematicians. In future works, we plan to prove these conjectures. As seen in the formulas presented in this paper, mathematical software such as Maple or Mathematica is essential for obtaining numerical results from our residue integrals. We plan to publish PDF copies of the Maple program worksheet on Masao Jinzenji's ResearchGate page \cite{JinR}.

\vspace{1em}

{\bf Acknowledgment} We would like to thank Prof. Y. Kanda, Prof. S. Kobayashi, Prof. H. Iritani and Prof. K. Ono for valuable discussions. Our research is partially supported by JSPS grant No. 22K03289. Research of M.J. is partially supported by JSPS grant No. 24H00182.

\section{Setup of Our Computation}
Let $M_{N}^{k}$ be a non-singular degree $k$ hypersurface in $CP^{N-1}$. Its complex dimension is clearly $N-2$. In this paper, we compute genus $1$ (elliptic) Gromov-Witten invariants of 
$M_{N}^{k}$ for $N\geq k$ by using the mirror map constructed from the multi-point virtual structure constants of $M_{N}^{k}$ \cite{JS1}.  
We denote by $h\in H^{1,1}(M_{N}^{k},{\bf C})$ the hyperplane class of $M_{N}^{k}$. 

Let us briefly review the conjectures proposed in \cite{JS1}. We first introduce the following polynomials:
\ba
e_{k}(x,y)&:=&\prod_{j=0}^{k}\left(jx+(k-j)y\right),\no\\
w_{a}(x,y)&:=&\frac{x^{a}-y^{a}}{x-y}=\sum_{j=0}^{a-1}x^jy^{a-1-j},\no\\
\ea
which serve as the  building blocks of the residue integral representation of genus $0$ multi-point  virtual structure constant 
$w({\cal O}_{h^{a}}{\cal O}_{h^{b}}|\prod_{p=0}^{N-2}({\cal O}_{h^{p}})^{m_{p}})_{0,d}$ of degree $d$.  For $d\geq 1$, it is explicitly given as follows.
\ba
&&w({\cal O}_{h^{a}}{\cal O}_{h^{b}}|\prod_{p=0}^{N-2}({\cal O}_{h^{p}})^{m_{p}})_{0,d}\no\\
&=&\frac{1}{(2\pi\sqrt{-1})^{d+1}}\oint_{C_{z_{0}}}dz_{0}\oint_{C_{z_{1}}}dz_{1}\cdots\oint_{C_{z_{d}}}dz_{d}\no\\
&&\times (z_{0})^{a}\left(\frac{\prod_{j=1}^{d}e_{k}(z_{j-1},z_{j})}{\prod_{i=1}^{d-1}kz_{i}(2z_{i}-z_{i-1}-z_{i+1})}\right)(z_{d})^{b}
\left(\prod_{p=0}^{N-2}\left(\sum_{n=1}^{d}w_{p}(z_{n-1},z_{n})\right)^{m_{p}}\right)                         \no\\
&&\times\prod_{q=0}^{d}\frac{1}{(z_{q})^{N}}.
\label{multivir}
\ea
In the residue integral above, $\oint_{C_{z_{0}}}dz_{0}$ and $\oint_{C_{z_{d}}}dz_{d}$ denote the  operation of taking residues at $z_{0}=0$ and $z_{d}=0$, respectively.
On the other hand, $\oint_{C_{z_{j}}}dz_{j}\;\;(j=1,2,\cdots,d-1)$ indicates taking residues at $z_{j}=0$ and $z_{j}=\frac{z_{j-1}+z_{j+1}}{2}$.  In this context, the result of the integration 
does not depend on the order in which the integrations are performed. If we take the residues in ascending order of the subscripts, the procedure is as follows. First, we compute the residue of the integrand at $z_{0}=0$, then add the residues of the resulting integrand at $z_{1}=0$ and $z_{1}=\frac{z_{2}}{2}$. Next, we add the residues of the subsequent integrand at $z_{2}=0,\frac{z_{3}}{2}$, and $z_{2}=\frac{2z_{3}}{3}$, and so forth. We recommend that readers consult our Maple files in \cite{JinR} for further details.    

If $d=0$, the integral is $0$ except for the following cases:
\ba
w({\cal O}_{h^{a}}{\cal O}_{h^{b}}|{\cal O}_{h^{c}})_{0,0}=k\delta_{a+b+c,N-2}.
\label{zero}
\ea
For $d\geq 1$, the following condition holds:
\ba
w({\cal O}_{h^{a}}{\cal O}_{h^{b}}|\prod_{p=0}^{N-2}({\cal O}_{h^{p}})^{m_{p}})_{0,d}\neq 0 \Longrightarrow a+b+\sum_{p=0}^{N-2}(p-1)m_{p}=N-3+(N-k)d.
\ea

\begin{figure}[H]
 \begin{center}
   \includegraphics[width=100mm]{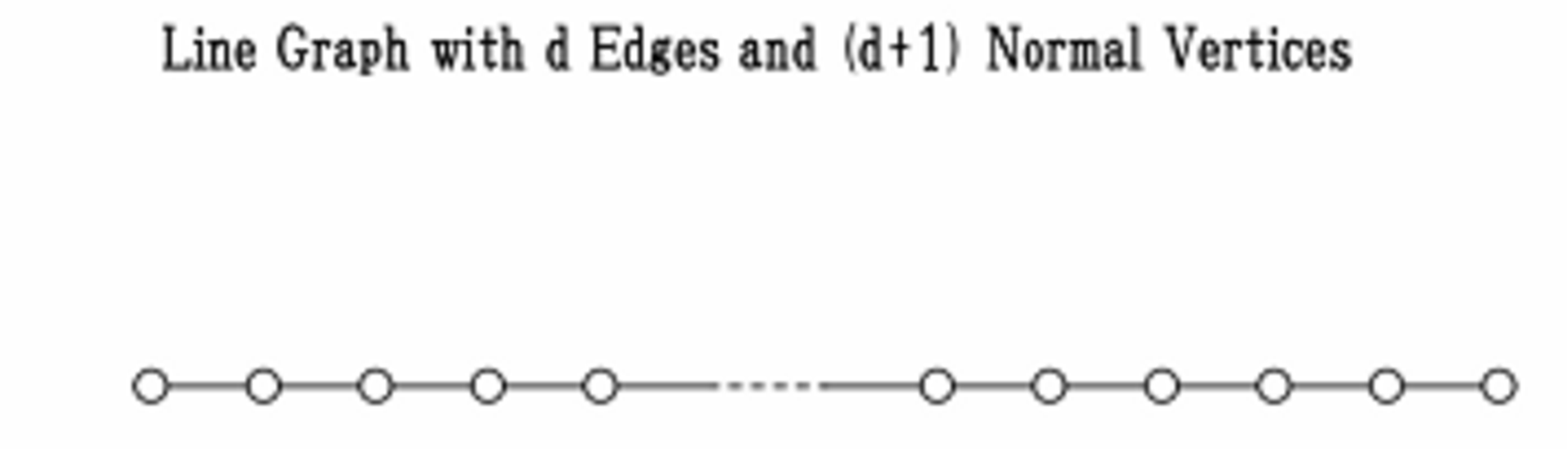}
\end{center}
  \label{line}
  \caption{The Graph used in Computing Genus $0$ Degree $d$ Virtual Structure Constants}
\end{figure}

We can easily see the following characteristics of  $w({\cal O}_{h^{a}}{\cal O}_{h^{b}}|\prod_{p=0}^{N-2}({\cal O}_{h^{p}})^{m_{p}})_{0,d}$:
\ba
w({\cal O}_{h^{a}}{\cal O}_{h^{b}}|\prod_{p=0}^{N-2}({\cal O}_{h^{p}})^{m_{p}})_{0,d}=\delta_{m_{0},0}\cdot d^{m_{1}}w({\cal O}_{h^{a}}{\cal O}_{h^{b}}|\prod_{p=2}^{N-2}({\cal O}_{h^{p}})^{m_{p}})_{0,d}
\;\;(d\geq 1 ),
\label{pk}
\ea
since we have $w_{0}(x,y)=0$ and $w_{1}(x,y)=1$. Next, we introduce ``perturbed two-point function'':
\ba
w({\cal O}_{h^a}{\cal O}_{h^b})_{0}(x^{0},x^{1},\cdots,x^{N-2}):=
\sum_{m_{0}=0}^{\infty}\cdots\sum_{m_{N-2}=0}^{\infty}\sum_{d=0}^{\infty}w({\cal O}_{h^{a}}{\cal O}_{h^{b}}|\prod_{p=0}^{N-2}({\cal O}_{h^{p}})^{m_{p}})_{0,d}\prod_{q=0}^{N-2}\frac{(x^{q})^{m_q}}{m_{q}!}.
\label{pert1}
\ea
(\ref{zero}) and (\ref{pk}) lead us to the following simplification:
\ba
&&w({\cal O}_{h^a}{\cal O}_{h^b})_{0}(x^{0},x^{1},\cdots,x^{N-2})\no\\
&=&kx^{N-2-a-b}+\sum_{m_{2}=0}^{\infty}\cdots\sum_{m_{N-2}=0}^{\infty}\sum_{d=1}^{\infty}e^{dx^{1}}w({\cal O}_{h^{a}}{\cal O}_{h^{b}}|\prod_{p=2}^{N-2}({\cal O}_{h^{p}})^{m_{p}})_{0,d}\prod_{q=2}^{N-2}\frac{(x^{q})^{m_q}}{m_{q}!}.
\label{simp}
\ea
Therefore, we only consider $w({\cal O}_{h^{a}}{\cal O}_{h^{b}}|\prod_{p=2}^{N-2}({\cal O}_{h^{p}})^{m_{p}})_{0,d}$ from now on.
In \cite{JS1}, we defined ``mirror map'', or coordinate change of deformation variables as follows:
\ba
t^p(x^{0},x^{1},\cdots,x^{N-2}):=\frac{1}{k}w({\cal O}_{h^{N-2-p}}{\cal O}_{1})_{0}\;\;\;(p=0,1,\cdots, N-2).
\ea 
(\ref{simp}) tells us that the above mirror map has the following structure:
\ba
t^p(x^{0},x^{1},\cdots,x^{N-2})=
x^{p}+\frac{1}{k}\sum_{m_{2}=0}^{\infty}\cdots\sum_{m_{N-2}=0}^{\infty}\sum_{d=1}^{\infty}e^{dx^{1}}w({\cal O}_{h^{N-2-p}}{\cal O}_{1}|\prod_{p=2}^{N-2}({\cal O}_{h^{p}})^{m_{p}})_{0,d}\prod_{q=2}^{N-2}\frac{(x^{q})^{m_q}}{m_{q}!}.
\ea
This structure allows us to invert the mirror map:
\ba
x^{p}=x^{p}(t^{0},t^{1},\cdots,t^{N-2})\;\;\;(p=0,1,\cdots,N-2).
\label{invert}
\ea

At this stage, we introduce the Gromov-Witten invariant $\langle\prod_{a=0}^{N-2}({\cal O}_{h^a})^{m_{a}}\rangle_{g,d}$ of genus $g$ and degree $d$.
In this paper, the genus $g$ is restricted to $0$ or $1$. We omit a rigorous definition of this invariant here but mention some of its characteristics.
In general, it is non-zero only if the following condition is satisfied: 
\ba
&&\langle\prod_{a=0}^{N-2}({\cal O}_{h^a})^{m_{a}}\rangle_{0,d}\neq 0\;\Longrightarrow \sum_{a=0}^{N-2}m_{a}(a-1)=N-5+d(N-k),\;\;\no\\
&&\langle\prod_{a=0}^{N-2}({\cal O}_{h^a})^{m_{a}}\rangle_{1,d}\neq 0\;\Longrightarrow \sum_{a=0}^{N-2}m_{a}(a-1)=d(N-k).\;\;\no\\
\label{sel}
\ea  
If $d=0$ and $g=0$, it is zero except for the following case:
\ba
\langle{\cal O}_{h^a}{\cal O}_{h^{b}}{\cal O}_{h^{c}}\rangle_{0,d}=k\cdot\delta_{N-2,a+b+c}.
\ea
If $d=0$ and $g=1$, it is zero except for the following case:
\ba
\langle{\cal O}_{h}\rangle_{1,0}=-\frac{1}{24}\int_{M_{N}^{k}}h\wedge c_{N-3}(T^{\prime}M_{N}^{k}),
\ea
where $c_{N-3}(T^{\prime}M_{N}^{k})$ is the second-top Chern class of $T^{\prime}M_{N}^{k}$, the homorphic tangent bundle of $M_{N}^{k}$. 
If $d\geq 1$, it satisfies the following equality:
\ba
\langle\prod_{a=0}^{N-2}({\cal O}_{h^a})^{m_{a}}\rangle_{g,d}=\delta_{m_0,0}d^{m_{1}}\langle\prod_{a=2}^{N-2}({\cal O}_{h^a})^{m_{a}}\rangle_{g,d}.
\ea
We restate the conjecture on genus $0$ Gromov-Witten invariants of $M_{N}^{k}$ proposed in \cite{JS1}.
As in (\ref{pert1}), we introduce the perturbed two-point genus $0$ Gromov-Witten invariant:
\ba
&&\langle{\cal O}_{h^a}{\cal O}_{h^b}\rangle_{0}(t^{0},t^{1},\cdots,t^{N-2})\no\\
&:=&\sum_{m_{0}=0}^{\infty}\cdots\sum_{m_{N-2}=0}^{\infty}\sum_{d=0}^{\infty}\langle{\cal O}_{h^{a}}{\cal O}_{h^{b}}\prod_{p=0}^{N-2}({\cal O}_{h^{p}})^{m_{p}}\rangle_{0,d}
\prod_{q=0}^{N-2}\frac{(t^{q})^{m_q}}{m_{q}!}\no\\
&=&kt^{N-2-a-b}+\sum_{m_{2}=0}^{\infty}\cdots\sum_{m_{N-2}=0}^{\infty}\sum_{d=1}^{\infty}e^{dt^{1}}\langle{\cal O}_{h^{a}}{\cal O}_{h^{b}}\prod_{p=2}^{N-2}({\cal O}_{h^{p}})^{m_{p}}\rangle_{0,d}\prod_{q=2}^{N-2}\frac{(t^{q})^{m_q}}{m_{q}!}.
\ea
Then the conjecture is stated as follows:
\begin{conj}{\bf \cite{JS1}}
\ba
\langle{\cal O}_{h^a}{\cal O}_{h^b}\rangle_{0}(t^{0},t^{1},\cdots,t^{N-2})=w({\cal O}_{h^a}{\cal O}_{h^b})_{0}(x^{0}(t^{*}),x^{1}(t^{*}),\cdots,x^{N-2}(t^{*})),
\ea
where $x^{p}(t^{*})$ is an abbreviation for $x^{p}(t^{0},t^{1},\cdots,t^{N-2})$ given in (\ref{invert}). 
\end{conj}
This conjecture provides a method for computing genus $0$ Gromov-Witten invariants of $M_{N}^{k}$ and it has been confirmed numerically for low degrees in many examples. 

In this paper, we also introduce the elliptic multi-point virtual structure constant $w(\prod_{p=0}^{N-2}({\cal O}_{h^{p}})^{m_{p}})_{1,d}$ of degree $d$, whose definition for the $d\geq 1$ case 
will be given in the next section. In the $d=0$ case, it is $0$ except for the following case:
\ba
w({\cal O}_{h})_{1,0}:=\langle{\cal O}_{h}\rangle_{1,0}=-\frac{1}{24}\int_{M_{N}^{k}}h\wedge c_{N-3}(T^{\prime}M_{N}^{k}). 
\label{w10}
\ea 
For $d\geq 1$, we can show the following characteristics of $w(\prod_{p=0}^{N-2}({\cal O}_{h^{p}})^{m_{p}})_{1,d}$  using the definition provided in the next section:
\ba
w(\prod_{p=0}^{N-2}({\cal O}_{h^{p}})^{m_{p}})_{1,d}=\delta_{m_0,0}d^{m_{1}}w(\prod_{p=2}^{N-2}({\cal O}_{h^{p}})^{m_{p}})_{1,d}.
\label{w1pk}
\ea
We then introduce the generating function of $w(\prod_{p=0}^{N-2}({\cal O}_{h^{p}})^{m_{p}})_{1,d}$'s:
\ba
&&F_{1,vir.}^{N,k, B}(x^{0},x^{1},\cdots,x^{N-2})\no\\
&:=&\sum_{d=0}^{\infty}\sum_{m_0=0}^{\infty}\cdots\sum_{m_{N-2}=0}^{\infty}w(\prod_{p=0}^{N-2}({\cal O}_{h^{p}})^{m_{p}})_{1,d}\prod_{q=0}^{N-2}\frac{(x^{q})^{m_q}}{m_{q}!}\no\\
&=&-\frac{1}{24}\left( \int_{M_{N}^{k}}h\wedge c_{N-3}(T^{\prime}M_{N}^{k}) \right)x^{1}+\sum_{d=1}^{\infty}\sum_{m_2=0}^{\infty}\cdots\sum_{m_{N-2}=0}^{\infty}e^{dx^{1}}w(\prod_{p=2}^{N-2}({\cal O}_{h^{p}})^{m_{p}})_{1,d}\prod_{q=2}^{N-2}\frac{(x^{q})^{m_q}}{m_{q}!}.
\ea
In the step from the second line to the third line, we used (\ref{w10}) and (\ref{w1pk}). On the other hand, we also define the generating function of genus 1 Gromov-Witten invariants of 
$M_{N}^{k}$.
\ba
&&F_{1}^{N,k,A}(t^{0},t^{1},\cdots,t^{N-2})\no\\
&:=&\sum_{d=0}^{\infty}\sum_{m_0=0}^{\infty}\cdots\sum_{m_{N-2}=0}^{\infty}\langle\prod_{p=0}^{N-2}({\cal O}_{h^{p}})^{m_{p}}\rangle_{1,d}\prod_{q=0}^{N-2}\frac{(t^{q})^{m_q}}{m_{q}!}\no\\
&=&-\frac{1}{24}\left( \int_{M_{N}^{k}}h\wedge c_{N-3}(T^{\prime}M_{N}^{k}) \right)t^{1}+\sum_{d=1}^{\infty}\sum_{m_2=0}^{\infty}\cdots\sum_{m_{N-2}=0}^{\infty}e^{dt^{1}}\langle\prod_{p=2}^{N-2}({\cal O}_{h^{p}})^{m_{p}}\rangle_{1,d}\prod_{q=2}^{N-2}\frac{(t^{q})^{m_q}}{m_{q}!}.
\ea
With this setup, we state our main conjecture in this paper.
\begin{conj}{\bf( Main Conjecture)}
\ba
F_{1}^{N,k,A}(t^{0},t^{1},\cdots,t^{N-2})=F_{1,vir.}^{N,k,B}(x^{0}(t^{*}),x^{1}(t^{*}),\cdots,x^{N-2}(t^{*})),
\ea
where $x^{p}(t^{*})$ is the inversion of the mirror map given in (\ref{invert}). 
\label{main}
\end{conj} 
With the explicit definition of the elliptic multi-point virtual structure constants provided in the next section, this conjecture offers a method for computing genus $1$ Gromov-Witten invariants of $M_{N}^{k}$. In Section 4, we test this conjecture by comparing our predictions with known results for genus $1$ Gromov-Witten invariants of $M_{N}^{k}$. 

\section{Perturbative Definitions of Elliptic Multi-Point Virtual Structure Constants }
\subsection{Graphs}
To explain the structure of graphs used in our computation, we first introduce the concept of a partition of a positive integer $d$.
\ba
\sigma=(d_{1},d_{2},\cdots,d_{l(\sigma)})\;\;\;(d_{1}\geq d_{2}\geq \cdots d_{l(\sigma)}>0,\;\;\sum_{i=1}^{l(\sigma)}d_{i}=d),
\ea
where we call $l(\sigma)$ the length of the partition $\sigma$. We denote by $P_{d}$ the set of all partitions of the positive integer $d$.
\ba
P_{d}:=\{\sigma=(d_{1},d_{2},\cdots,d_{l(\sigma)})\;|\;d_{1}\geq d_{2}\geq \cdots d_{l(\sigma)}>0,\;\;\sum_{i=1}^{l(\sigma)}d_{i}=d\;\}.
\ea
For later use, we define the symmetry factor associated with $\sigma\in P_{d}$:
\ba
\mbox{Sym}(\sigma):=\frac{(l(\sigma)-1)!}{\prod_{j=1}^{d}\mbox{mul}(\sigma;j)!}\;\;\;(\sigma\in P_{d}),
\ea
where $\mbox{mul}(\sigma;j)$ is the multiplicity of $i\;\;(1\leq i\leq d)$  in $\sigma=(d_{1},d_{2},\cdots,d_{l(\sigma)})$.
We also define the following rational number:
\ba
R_{N,k}(d)&:=&\frac{(N-1)}{2d}-\left(N-\frac{1}{k}\right)\frac{1}{d^2},
\ea
which plays the role of a symmetry factor for specific graphs.

In our computation, we use graphs with normal edges represented by ``---'' and three types of vertices.
\begin{itemize}
\item[(i)] {Normal Vertex}
\item[(ii)] {Elliptic Vertex}
\item[(iii)] {Cluster Vertex of degree $d$}
\end{itemize}
These are graphically represented by the following symbols.

\begin{figure}[H]
 \begin{center}
   \includegraphics[width=40mm]{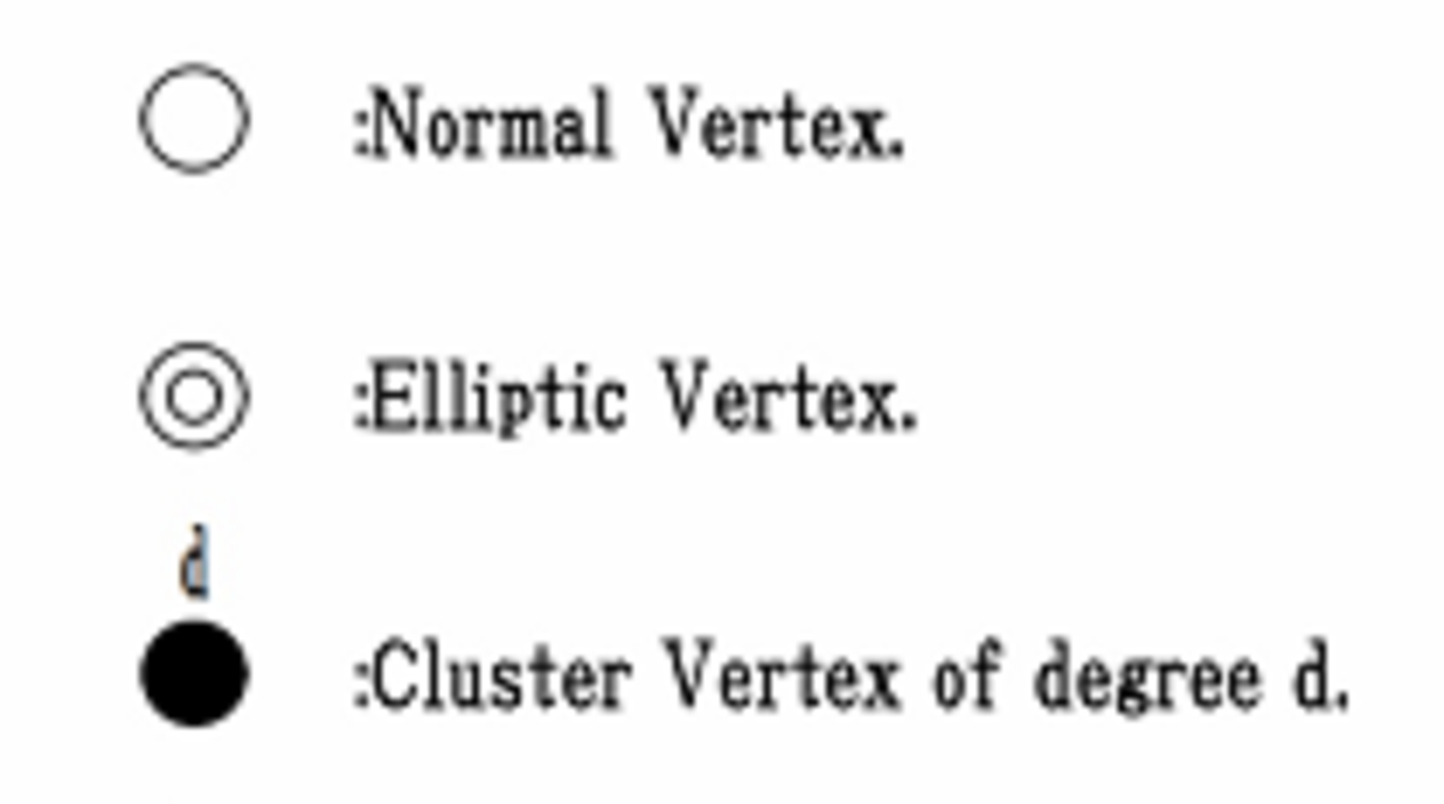}
\end{center}
  \label{ver}
  \caption{Vertices used in Our Construction}
\end{figure}

A single edge is assigned a degree of $1$. A single normal vertex, an elliptic
vertex and a cluster vertex of degree $d$ are assigned  degrees of $0$, $0$ and $d$, respectively. The graphs used in computing elliptic virtual constants of degree $d$ are classified into the following four types.
\begin{itemize}
\item[(i)] {Star graph associated with $\sigma\in P_{d}$ having an elliptic vertex as its center}  
\item[(ii)] {Loop graph with $d$ edges and $d$ normal vertices $(d\geq 2)$}
\item[(iii)] {Star graph associated with $\sigma\in P_{d-f}\;\;(1\leq f\leq d-1)$ having a cluster vertex of degree $f$ as its center} 
\item[(iv)] {Graph  consisting of a single cluster vertex of degree $d$}
\end{itemize}
The type (i) and type (ii) graphs appear as the "linearization" of graphs used in localization computations of genus 1 Gromov-Witten invariants of projective space \cite{GP,KZ}. The type (iii) and type (iv) graphs are introduced to cancel the "diagonal anomaly" mentioned in Section 1, which arises when rewriting the results of localization computations into the form of residue integrals. Examples of these four types are shown in the following figures.

\begin{figure}[H]
\begin{center}
   \includegraphics[width=52mm]{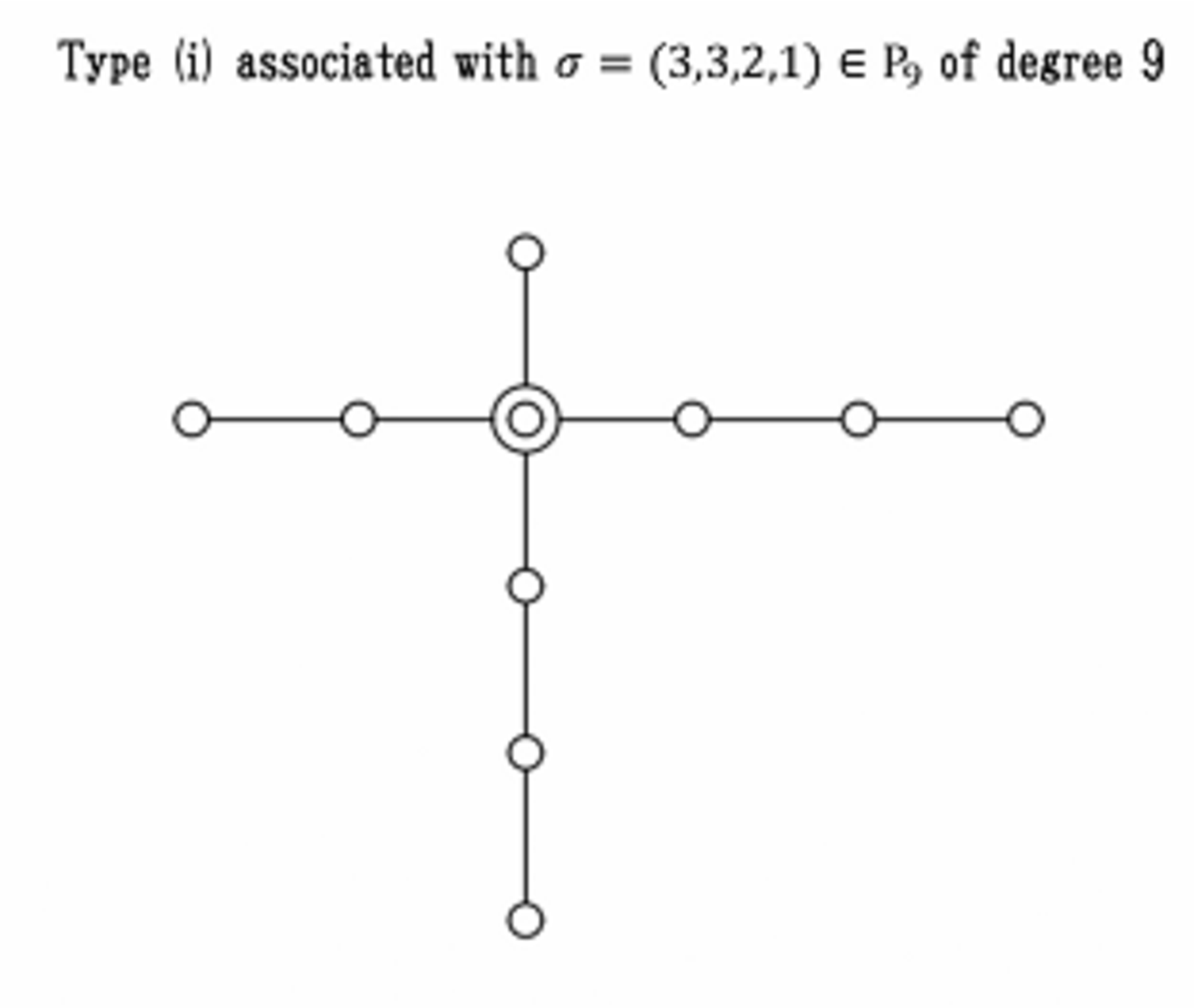}
\end{center}  
  \label{(i)}
  \caption{An example of Type (i) Graph }
\end{figure}

\begin{figure}[H]
 \begin{center}
   \includegraphics[width=55mm]{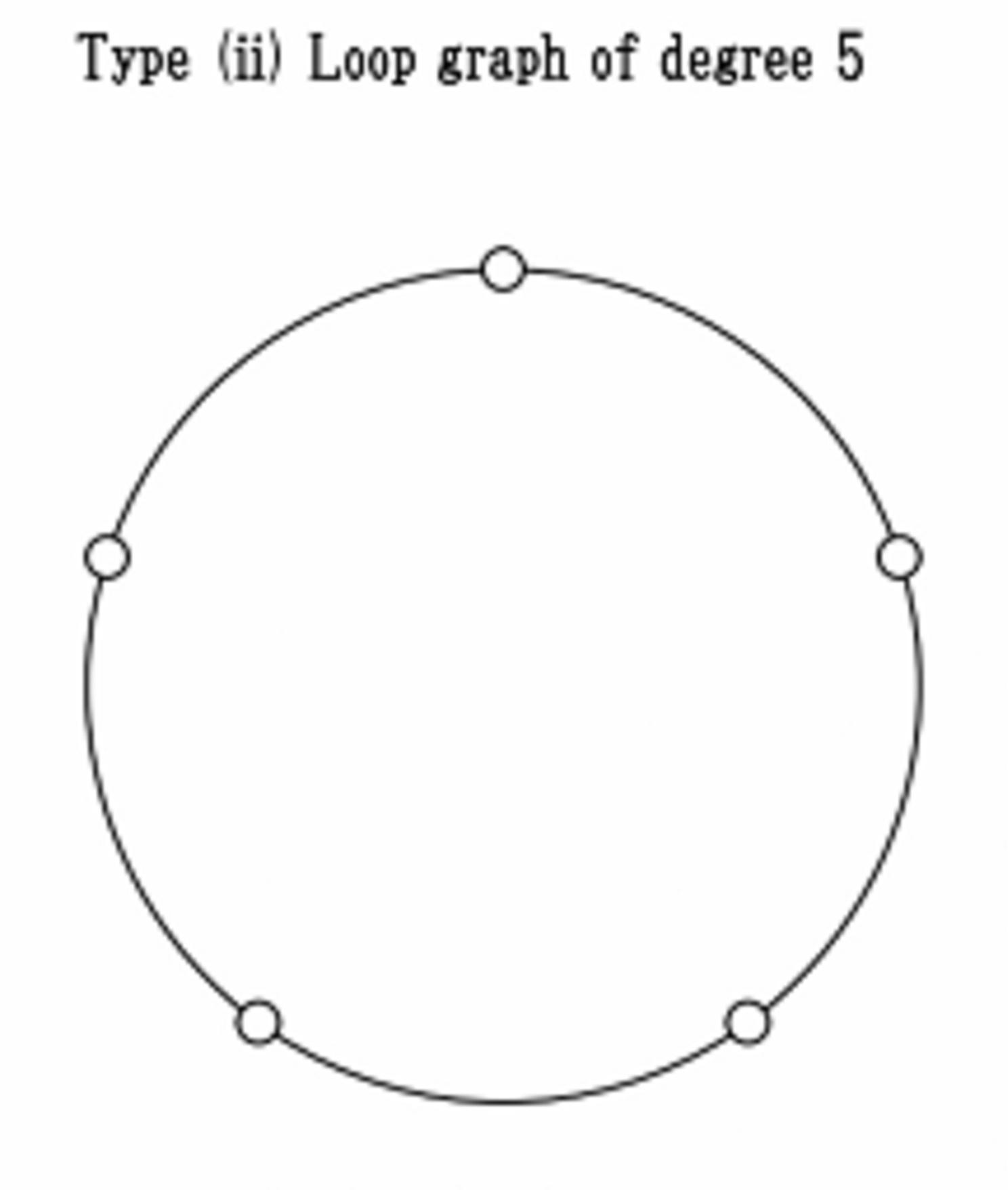}
\end{center}
  \label{(ii)}
  \caption{An example of Type (ii) Graph }
\end{figure}

\begin{figure}[H]
 \begin{center}
   \includegraphics[width=85mm]{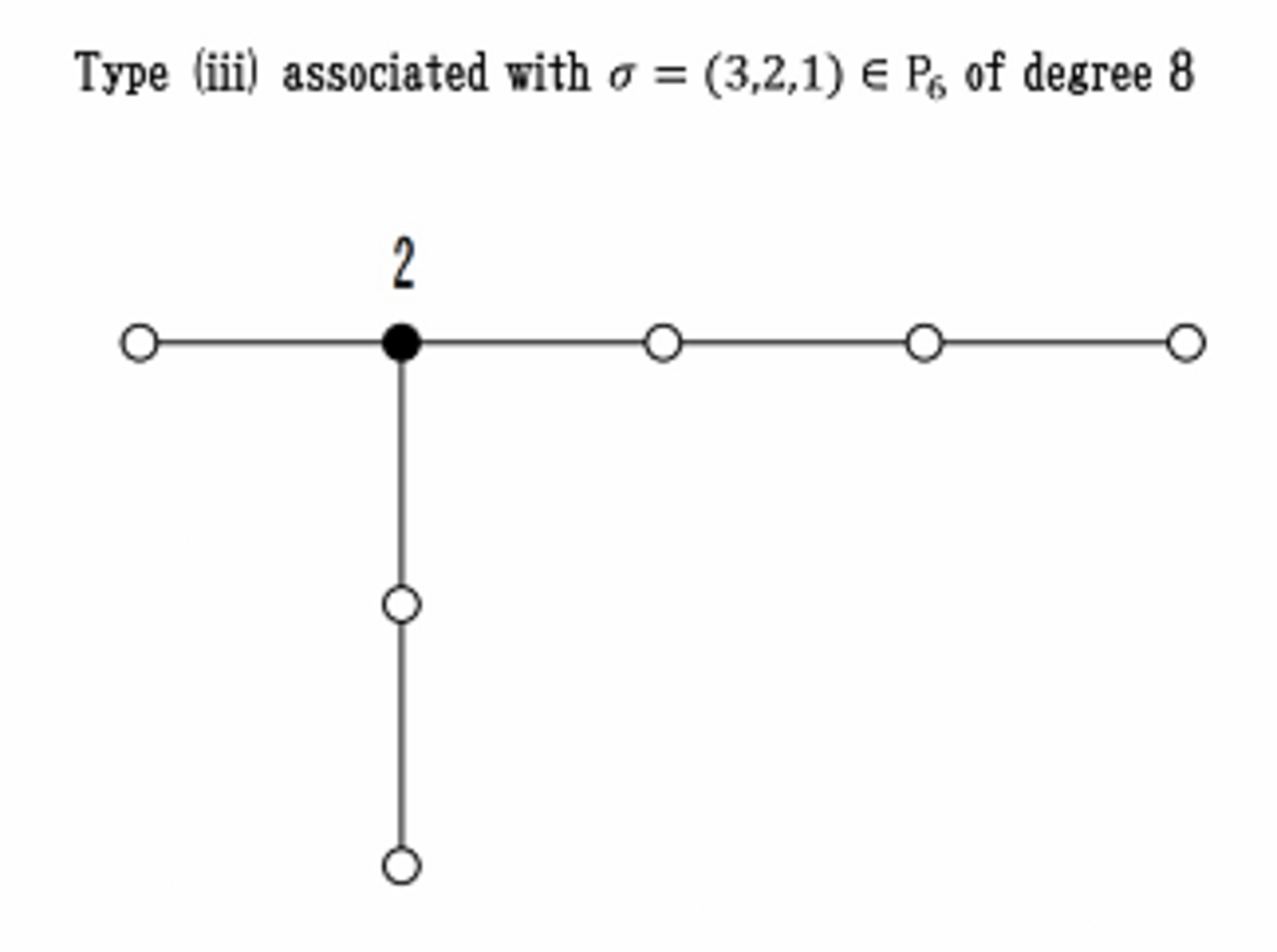}
\end{center}
  \label{(iii)}
  \caption{An example of Type (iii) Graph }
\end{figure}

\begin{figure}[H]
 \begin{center}
   \includegraphics[width=40mm]{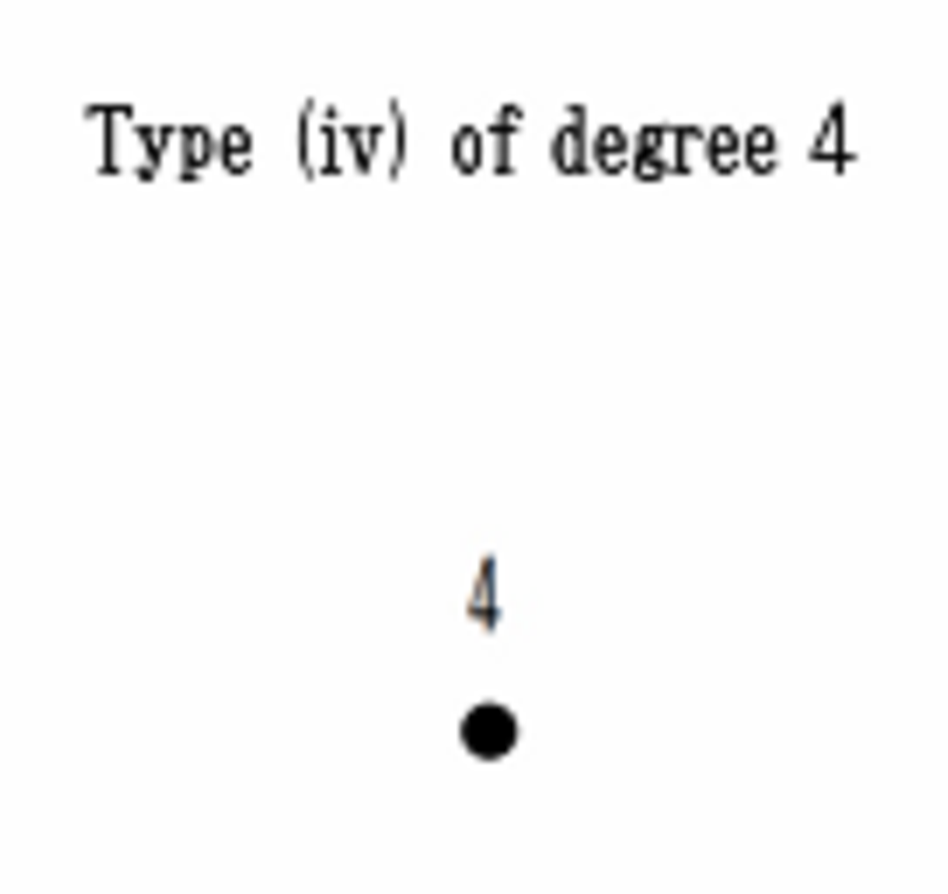}
\end{center}
  \label{(iv)}
  \caption{An example of Type (iv) Graph }
\end{figure}

\vspace{5cm}

We then list all the graphs used in computing elliptic virtual constants for $d=1,2,3$.

\begin{figure}[H]
 \begin{center}
   \includegraphics[width=70mm]{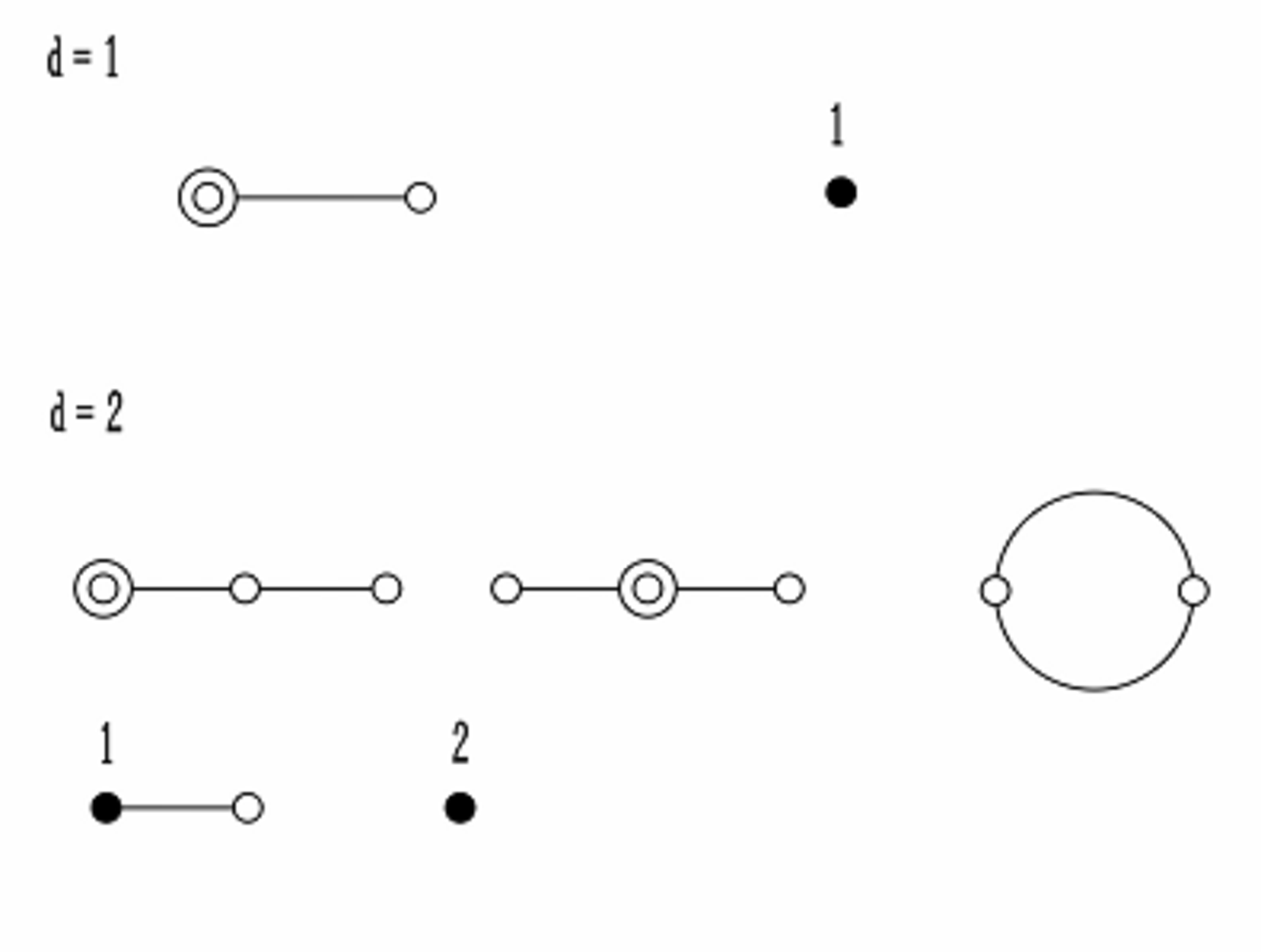}
\end{center}
  \label{d=1,2}
  \caption{Graphs used in Computing Elliptic Degree $1,2$ Virtual Structure Constants}
\end{figure}

\begin{figure}[H]
 \begin{center}
   \includegraphics[width=90mm]{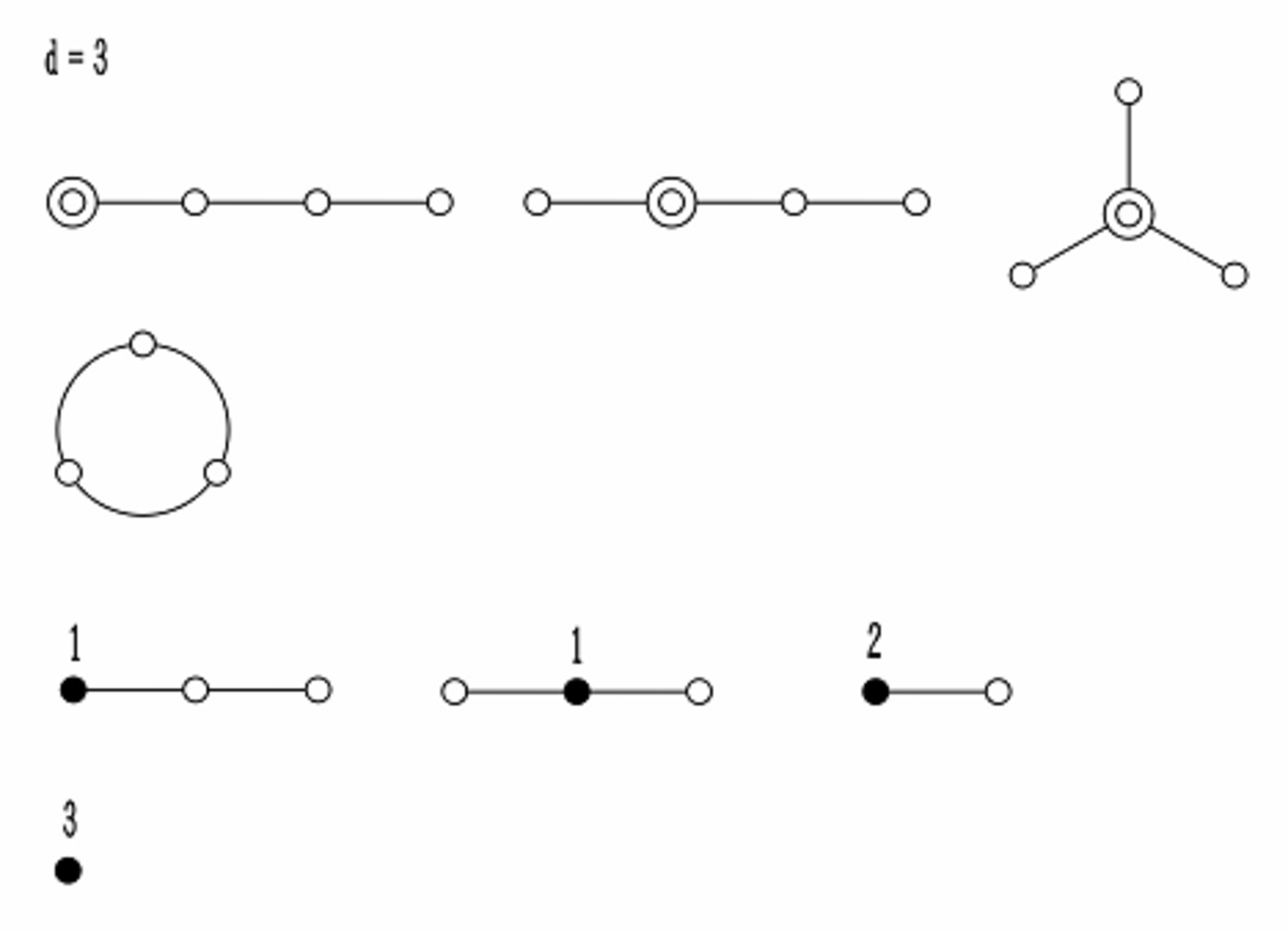}
\end{center}
  \label{d=3}
  \caption{Graphs used in Computing Elliptic Degree $3$ Virtual Structure Constants}
\end{figure}

We can easily see that the number of graphs used in computing elliptic virtual constants of degree $d$ (which we denote by $N_{d}$)is given by:
\ba
N_{d}=\left\{\begin{array}{cc}2&(d=1),\\ \sum_{j=1}^{d}(P_{j})^{\sharp} +1&(d\geq 2).\end{array}\right.
\ea
Therefore, the generating function of $N_{d}$ is given by:
\ba
1+\sum_{d=1}^{\infty}N_{d}q^{d}&=&\frac{q^2}{1-q}+\frac{1}{1-q}\left(\prod_{m=1}^{\infty}\frac{1}{1-q^{m}}\right)\no\\
&=&1+2q+5q^2+8q^3+13q^4+20q^5+31q^6+46q^7+68q^8+98q^9+140q^{10}+\cdots.
\ea
We can see that the number of graphs increases rather slowly as $d$ rises. 

\subsection{Integrands of Residue Integrals}

Elliptic virtual structure constant $w(\prod_{a=0}^{N-2}({\cal O}_{h^a})^{m_{a}})_{1,d}$ is computed by summing up contributions from these graphs of degree $d$. 
Contribution from each graph is given as residue integral of multi-variables determined from the graph. Let us write down the integrand of the residue integral associated with 
the four types of graphs.

We begin with type (i) graph associated with $\sigma=(d_{1},d_{2},\cdots, d_{l})\in P_{d}$. In this case, we prepare $(d+1)$ complex variables $z_{0}$ and $z_{i,j}\;\;(1\leq i\leq l,\;\;1\leq j\leq 
d_{i})$. $z_{0}$ is associated with the center elliptic vertex, and $z_{i,j}$ is associated with a vertex in the tail of a star graph with $d_{i}$ edges. It is well-known that the total Chern class of $T^{\prime}M_{N}^{k}$ is given by 
\ba
c(T^{\prime}M_{N}^{k})=\frac{(1+h)^N}{1+kh},
\ea
and the top Chern class $c_{N-2}(T^{\prime}M_{N}^{k})$ is given as  $h^{N-2}$ multiplied by coefficient of $h^{N-2}$ in $c(T^{\prime}M_{N}^{k})$. Then we define $c_{T}(z)$, which is a monomial in $z$, by the following equality.
\ba
c_{T}(h)=c_{N-2}(T^{\prime}M_{N}^{k}).
\ea   
Then integrand associated with the graph is given as follows.
\ba
&&\mbox{Sym}(\sigma)\frac{1}{24}\frac{c_{T}(z_{0})}{(kz_{0})^{l-1}}\left(\prod_{i=1}^{l}\frac{e_{k}(z_{0},z_{i,1})}{z_{i,1}-z_{0}}\right)
\left(\prod_{i=1}^{l}\left(\prod_{j=1}^{d_{i}-1}\frac{e_{k}(z_{i,j},z_{i,j+1})}{(2z_{i,j}-z_{i,j-1}-z_{i,j+1})kz_{i,j}}\right)\right)\no\\
&&\times \left(\frac{1}{(z_{0})^{N}}\right)\left(\prod_{i=1}^{l}\prod_{j=1}^{d_{i}}\frac{1}{(z_{i,j})^N}\right)
\left(\prod_{a=2}^{N-2}\left(\sum_{i=1}^{l}\sum_{j=1}^{d_{i}}w_{a}(z_{i,j-1},z_{i,j})\right)^{m_a}\right),\no\\
&&\hspace{9cm}(z_{i,0}=z_{0}\; (i=1,2,\cdots,l)).
\ea
Integrand for type (ii) loop graph of degree $d$ is given by:
\ba
&&\frac{1}{2d}\left(\prod_{j=1}^{d}\frac{e_{k}(z_{j},z_{j+1})}{(2z_{j}-z_{j-1}-z_{j+1})kz_{j}}\right)\left(\prod_{a=2}^{N-2}\left(\sum_{j=1}^{d}w_{a}(z_{j},z_{j+1})\right)^{m_a}\right)\cdot
\left(\prod_{j=1}^{d}\frac{1}{(z_{j})^{N}}\right),\no\\
&&\hspace{9cm} (z_{0}=z_{d},\;z_{d+1}=z_{1}).
\ea
Integrand for  type (iii) graph of degree $d$ with a cluster vertex of degree $f$ $(1\leq f\leq d-1)$ associated with $\sigma=(d_{1},d_{2},\cdots, d_{l})\in P_{d-f}$ is given as
\ba
&&\mbox{Sym}(\sigma)\frac{1}{24}\left(-\frac{N-1}{N}\frac{1}{w^N}-\frac{N+1}{N}\frac{1}{(z_{0})^{N}}\right)\frac{1}{(w-z_{0})^{2}}\frac{e_{k}(w,z_{0})}{kw}
\left(\frac{e_{k}(z_{0},z_{0})}{kz_{0}}\right)^{f-1}\no\\
&&\times\frac{1}{(kz_{0})^{l-1}}\left(\prod_{i=1}^{l}\frac{e_{k}(z_{0},z_{i,1})}{z_{i,1}-z_{0}}\right)
\left(\prod_{i=1}^{l}\left(\prod_{j=1}^{d_{i}-1}\frac{e_{k}(z_{i,j},z_{i,j+1})}{(2z_{i,j}-z_{i,j-1}-z_{i,j+1})kz_{i,j}}\right)\right)\left(\frac{1}{(z_{0})^{N(f-1)}}\right)\left(\prod_{i=1}^{l}\prod_{j=1}^{d_{i}}\frac{1}{(z_{i,j})^N}\right)\no\\
&& \times\left(\prod_{a=2}^{N-2}\left(w_{a}(w,z_{0})+(f-1)w_{a}(z_{0},z_{0})+\sum_{i=1}^{l}\sum_{j=1}^{d_{i}}w_{a}(z_{i,j-1},z_{i,j})\right)^{m_a}\right),\no\\
&&\hspace{9cm}(z_{i,0}=z_{0}\; (i=1,2,\cdots,l)).
\ea
Lastly, integrand for type (iv) graph of degree $d$ is given by,
\ba
&&\frac{1}{24}R_{N,k}(d)\left(\frac{e_{k}(z_{0},z_{0})}{kz_{0}}\right)^{d}\left(\prod_{a=2}^{N-2}\left(dw_{a}(z_{0},z_{0})\right)^{m_a}\right)\frac{1}{(z_{0})^{Nd+1}}\no\\
&=&\frac{1}{24}R_{N,k}(d)\frac{k^{kd}\prod_{a=2}^{N-2}(da)^{m_{a}}}{z_{0}}.
\label{point}
\ea
Let $\Gamma$ be a graph of degree $d$ introduced so far, and $\mbox{Graph}_{d}$ be the set of all the graphs of degree $d$ from type (i) to (iv). Then we denote by $f_{\Gamma}$ 
the integrand defined above. In the following, we define a map $\mbox{Res}:f_{\Gamma}\rightarrow {\bf R}$, which corresponds to the operation of taking the residue of 
$f_{\Gamma}$. 
\begin{defi}
$\mbox{Res}:f_{\Gamma}\rightarrow {\bf R}$ is defined for each type of graph as follows.
\begin{itemize}
\item[(i)] {We first take the residue of $f_{\Gamma}$ at $z_{0}=0$. Next we take the residue of the resulting function at $z_{i,j}=0$ and $z_{i,j}=\frac{z_{i,j-1}+z_{i,j+1}}{2}$ and add them
up sequentially in ascending order of $j$ $(1\leq j\leq d_{i}-1)$. Lastly, we take the residue of the resulting function at $z_{i,d_{i}}=0$. The order among different $i$'s does not matter.}
\item[(ii)]{We take the residue of $f_{\Gamma}$ at $z_{j}=0$ and $z_{j}=\frac{z_{j-1}+z_{j+1}}{2}$ and add them up sequentially in ascending order of $j$ $(1\leq j\leq d)$.}
\item[(iii)]{We first take the residue of $f_{\Gamma}$ at $w=z_{0}$. Then the remaining process is the same as for type (i).}   
\item[(iv)]{We take the residue of $f_{\Gamma}$ at $z_{0}=0$.}
\end{itemize}
\end{defi}
Then the elliptic multi-point virtual structure constant is defined as follows:
\begin{defi}
\ba
w(\prod_{a=2}^{N-2}({\cal O}_{h^a})^{m_{a}})_{1,d}:=\sum_{\Gamma\in \mbox{Graph}_{d}}\mbox{Res}(f_{\Gamma})\;\hspace{3cm}(d\geq 1).
\ea
\label{defw1}
\end{defi}
In the next section, we explicitly compute $w(\prod_{a=2}^{N-2}({\cal O}_{h^a})^{m_{a}})_{1,d}$'s  using this definition for various $M_{N}^{k}$'s and test Conjecture \ref{main}
by comparing our predictions with known results.


\section{Numerical Tests by Various Examples}
\subsection{Tests for Fano Hypersurfaces}
First, we test Conjecture \ref {main}  using the examples of the surface $CP^{2}=M_{4}^{1}$, quadric surface $M_{4}^{2}$, and cubic surface $M_{4}^{3}$.  
Here, we present the explicit computation process for  the case of $CP^{2}$. In this case, the  generating function of the elliptic virtual structure constants $F_{1,vir.}^{4,1, B}$
has the following structure: 
\ba
&&F_{1,vir.}^{4,1, B}(x^{0},x^{1},x^{2})=-\frac{1}{8}x^{1}+\sum_{d=1}^{\infty}e^{dx^{1}}w(({\cal O}_{h^{2}})^{3d})_{1,d}\frac{(x^{2})^{3d}}{(3d)!},
\ea
where we used $c(T^{\prime}M_{4}^{1})=(1+h)^3=1+3h+3h^2$ and the condition:
\ba
w(({\cal O}_{h^{2}})^{a})_{1,d}\neq 0\Longrightarrow a=(N-k)d=3d.
\ea
We abbreviate $w(({\cal O}_{h^{2}})^{a})_{1,d}$  as $w_{a}$ and present the results of the computation using Definition \ref{defw1} in Table \ref{m41}. 
Therefore, $F_{1,vir.}^{4,1, B}$ is explicitly given by: 
\ba
&&F_{1,vir.}^{4,1, B}(x^{0},x^{1},x^{2})\no\\
&=&-\frac{1}{8} x^1-\frac{1}{16} e^{x^{1}} (x^2)^3-\frac{7}{80} e^{2x^1} (x^2)^6-\frac{77789}{362880} e^{3x^1}( x^2)^9-\frac{21344159}{31933440} e^{4x^{1}} (x^2)^{12}
-\frac{15774542951}{6604416000} e^{5x^{1}} (x^2)^{15}-\cdots\no\\
\ea 
On the other hand, the mirror map is evaluated using (\ref{multivir}) as follows. 
\ba
t^{0} &=& x^0 + \frac{1}{2}(x^2)^2e^{x^{1}} + \frac{8}{15}(x^2)^5e^{2x^{1}} + \frac{983}{840}(x^2)^8e^{3x^{1}} + \frac{4283071}{1247400}(x^2)^{11}e^{4x^{1}} + \frac{4019248213}{340540200}(x^2)^{14}e^{5x^{1}}+\cdots,\no\\
t^{1}&=& x^1+\frac{1}{2} (x^2)^3 e^{x^{1}}+\frac{7}{10} (x^2)^6 e^{2x^{1}}+\frac{2593}{1512} (x^2)^9 e^{3x^{1}}+\frac{2668063}{498960}(x^2)^{12} e^{4x^{1}}
+\frac{120501923}{6306300} (x^2)^{15} e^{5x^{1}}+\cdots,\no\\
t^{2}&=& x^2 + \frac{1}{4}(x^2)^4e^{x^{1}} + \frac{33}{70}(x^2)^7e^{2x^{1}} + \frac{16589}{12600}(x^2)^{10}e^{3x^{1}} + \frac{143698921}{32432400}(x^2)^{13}e^{4x^{1}} + 
\frac{75631936691}{4540536000}(x^2)^{16}e^{5x^{1}}+\cdots.\no\\
\label{mirr41}
\ea
In order to obtain $F_{1}^{4,1, A}(t^{0},t^{1},t^{2})$ from Conjecture \ref{main}, we  need $x^{1}=x^{1}(t^{1},t^{2})$ and $x^{2}=x^{2}(t^{1},t^{2})$ since both $F_{1,vir.}^{4,1, B}$ and 
$F_{1}^{4,1,A}$ do not depend on $x^{0}$ and $t^{0}$. After inverting (\ref{mirr41}), they are given as follows:
\ba
x^1&= &t^1-\frac{1}{2} (t^2)^3 e^{t^{1}}-\frac{3}{40} (t^2)^6 e^{2t^{1}}-\frac{3827}{30240} (t^2)^9 e^{3t^{1}}-\frac{4914517}{19958400}(t^2)^{12} e^{4t^{1}}
-\frac{10460869973}{18162144000}(t^2)^{15} e^{5t^{1}}-\cdots\no\\
x^2 &= &t^2-\frac{1}{4} (t^2)^4 e^{t^{1}}-\frac{27}{280}(t^2)^7 e^{2t^{1}}-\frac{7811}{50400}(t^2)^{10} e^{3t^{1}}-\frac{82505777}{259459200} (t^2)^{13} e^{4t^{1}}-
\frac{1014012107}{1320883200} (t^2)^{16} e^{5t^{1}}-\cdots,\no\\
\label{inv41}
\ea 
By substituting (\ref{inv41}) into $F_{1,vir.}^{4,1, B}(x^{1},x^{2})$, we obtain from Conjecture 2,
\ba
F_{1}^{4,1,A}(t^{1},t^{2})&=&F_{1,vir.}^{4,1, B}(x^{1}(t^{1},t^{2}),x^{2}(t^{1},t^{2}))\no\\
&=&-\frac{1}{8}t^1 + \frac{1}{362880}(t^2)^9e^{3t^{1}} + \frac{1}{2128896}(t^2)^{12}e^{4t^{1}} + \frac{173}{2594592000}(t^2)^{15}e^{5t^{1}}+\cdots.
\ea 
Then we have,
\ba
&&\langle({\cal O}_{h^{2}})^{3}\rangle_{1,1}=0,\;\langle({\cal O}_{h^{2}})^{6}\rangle_{1,2}=0,\;\langle({\cal O}_{h^{2}})^{9}\rangle_{1,3}=\frac{1}{362880}\cdot 9!=1,\;
\;\langle({\cal O}_{h^{2}})^{12}\rangle_{1,4}=\frac{1}{2128896}\cdot 12!=225,\;\no\\
&&\langle({\cal O}_{h^{2}})^{15}\rangle_{1,5}=\frac{173}{2594592000}\cdot 15!=87192.
\ea
These results are summarized in Table \ref{m41}, where $\langle({\cal O}_{h^{2}})^{a}\rangle_{1,d}$ is abbreviated to $N^{1}_{d,a}$.
\begin{table}[H]
\centering
\caption{$M_4^1$}
\begin{tabular}{|c|l|l|l|}
\hline
d&a& $N^1_{d,a}$ &$w_a$  \\
\hline
1 & 3 &0&$-\frac{3}{8}$\\ \hline
2& 6 & 0&$-63$\\ \hline
3& 9& 1 &$-77789$\\ \hline
4& 12 & 225&$-320162385$\\ \hline
5& 15& 87192&$-3123359504298$\\ \hline
\end{tabular}
\label{m41}
\end{table}
These results coincide with those presented in the paper \cite{G} by Getzler. We computed up to $d=5$ due to the increasing computational effort required for evaluating residue integrals as the degree increases.

As for the quadric surface and cubic surface, we  present only structure of $F_{1,vir.}^{N,k, B}$'s, the mirror maps, and tables of $w_{a}$'s and $N_{d,a}^{1}$'s. 
The processes of computation are the same as those for $M_{4}^{1}$.

$F_{1,vir.}^{4,2, B}$ and $F_{1,vir.}^{4,3, B}$ have the following structure.
\ba
&&F_{1,vir.}^{4,2, B}(x^{0},x^{1},x^{2})=-\frac{1}{6}x^{1}+\sum_{d=1}^{\infty}e^{dx^{1}}w_{2d}\frac{(x^{2})^{2d}}{(2d)!},\no\\
&&F_{1,vir.}^{4,3, B}(x^{0},x^{1},x^{2})=-\frac{1}{8}x^{1}+\sum_{d=1}^{\infty}e^{dx^{1}}w_{d}\frac{(x^{2})^{d}}{d!}.
\ea
Mirror map of $M_{4}^{2}$ is given by:
\ba
t^{0} &=& x^0 + 2x^2e^{x^{1}} + 10(x^2)^3e^{2x^{1}} + \frac{320}{3}(x^2)^5e^{3x^{1}} + \frac{53856}{35}(x^2)^7e^{4x^{1}} + \frac{74056288}{2835}(x^2)^9e^{5x^{1}}+\cdots,\no\\
t^1 &=& x^1 + 3(x^2)^2e^{x^{1}} + \frac{131}{6}(x^2)^4e^{2x^{1}} + \frac{12329}{45}(x^2)^6e^{3x^{1}} + \frac{121475}{28}(x^2)^8e^{4x^{1}} +\frac{370005883}{4725}(x^2)^{10}e^{5x^{1}}+
\cdots,\no\\
t^{2}&=&x^2 + 2(x^2)^3e^{x^{1}} + \frac{313}{15}(x^2)^5e^{2x^{1}} + \frac{10764}{35}(x^2)^7e^{3x^{1}} + \frac{15178391}{2835}(x^2)^9e^{4x^{1}} + \frac{458817242}{4455}(x^2)^{11}e^{5x^{1}}
+\cdots, \no\\
\ea
and that of  $M_{4}^{3}$ is given as follows.
\ba
t^0 &=&x^{0}+6e^{x^{1}}++ 144x^{2}e^{2x^{1}}++ 7398(x^2)^2e^{3x^{1}}+520344(x^2)^3e^{4x^{1}}+43392510(x^2)^4e^{5x^{1}} +\cdots,\no\\
t^{1} &=& x^1 + 21x^2e^{x^{1}} + \frac{1611}{2}(x~2)^2e^{2x^{1}} + 52191(x^2)^3e^{3x^{1}} + \frac{16915311}{4}(x^2)^4e^{4x^{1}} + 388313757(x^2)^5e^{5x^{1}}+\cdots,\no\\
t^2 &= &x^2 + 21(x^2)^2e^{x^{1}} + 1305(x^2)^3e^{2x^{1}} + 106056(x^2)^4e^{3x^{1}} + \frac{49255533}{5}(x^2)^5e^{4x^{1}} + \frac{4964539329}{5}(x^2)^6e^{5x^{1}}+\cdots.\no\\
\ea
Tables of  $w_{a}$'s and $N_{d,a}^{1}$'s are given as follows.
\begin{table}[H]
\centering
\caption{$M_4^2$}
\begin{tabular}{|c|l|l|l|l|}
\hline
d&a& $N^1_{d,a}$ &$N^1_{d,a}/2^{a}$&$w_a$  \\
\hline
1 & 2 &0&0&$-1$\\ \hline
2& 4 & 0& 0&$-\frac{262}{3}$ \\ \hline
3& 6& 0 & 0&$-\frac{98632}{3}$ \\ \hline
4& 8 & 256& 1&$-29153744$\\ \hline
5& 10& 40960 & 40 &$-47360066944$ \\ \hline
\end{tabular}
\end{table}

\begin{table}[H]
\centering
\caption{$M_4^3$}
\begin{tabular}{|c|l|l|l|l|}
\hline
d&a& $N^1_{d,a}$ &$N^1_{d,a}/3^{a}$ &$w_a$ \\
\hline
1 & 1 &0&0&$-\frac{21}{8}$\\ \hline
2& 1 & 0& 0 &$-\frac{1611}{8}$\\ \hline
3& 3& 27 & 1&$-\frac{156465}{4}$ \\ \hline
4& 4 &2187& 27&$-\frac{50682753}{4}$\\ \hline
5& 5 & 183708 & 756 &$-5815337247$ \\ \hline
\end{tabular}
\end{table}

In these tables, $N_{d,a}^{1}/k^{a}$ values are presented because Poincare dual of a point in $M_{N}^{k}$ is given by $\frac{1}{k}h^{2}$ instead of $h^{2}$. Therefore, $N_{d,a}^{1}/k^{a}$ values
turn out to be integers.

Next, we test Conjecture \ref{main} for $N=5$ cases, i.e., complex Fano 3-folds: $M_{5}^{1}, M_{5}^{2}, M_{5}^{3}, M_{5}^{4}$. In these cases, we compute the elliptic virtual structure 
constant $w(({\cal O}_{h^{2}})^{a}({\cal O}_{h^{3}})^{b})_{1,d}$ and the genus $1$ Gromov-Witten invariant $\langle ({\cal O}_{h^{2}})^{a}({\cal O}_{h^{3}})^{b}\rangle_{1,d}$. These are non-zero 
only if the following condition is satisfied:
\ba
a+2b=(N-k)d\;(=(5-k)d).
\label{sel3}
\ea
We abbreviate these as $w_{a,b}$ and $N_{d,a,b}^{1}$, and present results of numerical computations using Conjecture \ref{main} in the tables collected in Appendix A.  
For the 3-fold cases, $N_{d,a,b}^{1}$ is not always an integer due to degenerate contributions from genus $0$ Gromov-Witten invariants $N_{d,a,b}^{0}=
\langle ({\cal O}_{h^{2}})^{a}({\cal O}_{h^{3}})^{b}\rangle_{0,d}$. This contribution is non-zero only if (\ref{sel3}) is satisfied. According to the paper \cite{P} by Pandharipande, the contribution is 
given by: 
\ba
-\frac{1}{24}((N-k)d-2)N_{d,a,b}^{0}.
\ea
Hence we expect that the expression
\ba
\frac{1}{24}((N-k)d-2)N_{d,a,b}^{0}.+N_{d,a,b}^{1},
\ea
becomes an integer because it is expected to count the number of genus $1$ curves in $M_{N}^{k}$ that satisfy ``passing-through'' conditions imposed by the operator insertions. 
We observe that this combination is indeed integer-valued in the tables presented in Appendix A.   
The results of $M_{5}^{1}=CP^{3}$ coincide with those presented in the paper \cite{G} by Getzler, and the results for $M_{5}^{3}$ coincide with those obtained from Virasoro Conjecture 
by Eguchi et al. \cite{EJX}. In comparing our results with those  in \cite{EJX}, we must divide $N_{d,a,b}^{1}$ and $N_{d,a,b}^{0}$ by $k^{a+b}=3^{a+b}$. This is because 
fundamental integral cohomology classes of $H^{2,2}(M_{5}^{k},{\bf C})$ and $H^{3,3}(M_{5}^{k},{\bf C})$ are given by $\frac{1}{k}h^{2}$ and $\frac{1}{k}h^{3}$, respectively. In \cite{EJX},
$\frac{1}{k}h^{2}$ and $\frac{1}{k}h^{3}$ are used as cohomology class  corresponding to operator insertions, whereas we use $h^{2}$ and $h^{3}$ instead. Thus, we divide our 
results by  $k^{a+b}$. We do not include the mirror maps used in deriving our results for Fano 3-folds, which we leave to readers as exercises.

\subsection{Tests for Calabi-Yau Hypersurfaces and Relation to BCOV-Zinger formula} 
In this subsection, we consider Calabi-Yau hypersurface $M_{k}^{k}$ in $CP^{k-1}$.
We introduce the following power series in $e^{x}$, which was originally introduced in \cite{CJ}.
\ba
\tilde{L}^{k,k}_{m}(e^{x}):=1+\sum_{d=1}^{\infty}\frac{w({\cal O}_{h^{k-2-m}}{\cal O}_{h^{m-1}}|{\cal O}_{h})_{0,d}}{k}e^{dx}.
\ea
In this case, we also introduce the two-point virtual structure constant defined by:
\ba
&&w({\cal O}_{h^{a}}{\cal O}_{h^{b}})_{0,d}\no\\
&=&\frac{1}{(2\pi\sqrt{-1})^{d+1}}\oint_{C_{z_{0}}}dz_{0}\oint_{C_{z_{1}}}dz_{1}\cdots\oint_{C_{z_{d}}}dz_{d}(z_{0})^{a}\left(\frac{\prod_{j=1}^{d}e_{k}(z_{j-1},z_{j})}{\prod_{i=1}^{d-1}kz_{i}(2z_{i}-z_{i-1}-z_{i+1})}\right)(z_{d})^{b}\left(\prod_{j=0}^{d}\frac{1}{(z_{q})^{N}}\right),\no\\
\label{vir}
\ea
where the integration paths are the same as those used in defining $w({\cal O}_{h^{a}}{\cal O}_{h^{b}}|\prod_{p=0}^{N-2}({\cal O}_{h^{p}})^{m_{p}})_{0,d}$.
It satisfies the condition:
\ba
w({\cal O}_{h^{a}}{\cal O}_{h^{b}})_{0,d}\neq 0\Longrightarrow a+b=k-3.
\ea
Then, we define a two-point virtual structure constant perturbed by $x$:
\ba
w({\cal O}_{h^{k-2-m}}{\cal O}_{h^{m-1}})_{0}(x):=kx+\sum_{d=1}^{\infty}e^{dx}w({\cal O}_{h^{k-2-m}}{\cal O}_{h^{m-1}})_{0,d}.
\ea
As shown in \cite{Jin1, Jin3}, the mirror map $t=t(x)$ used in the mirror computation of Gromov-Witten invariants of $M_{k}^{k}$ is given by:
\ba
t(x)=\frac{1}{k}w({\cal O}_{h^{k-3}}{\cal O}_{1})_{0}(x).
\label{mirr}
\ea  
In \cite{Jin2}, we proved that the Gromov-Witten invariant $\langle{\cal O}_{h^{k-2-m}}{\cal O}_{h^{m-1}}\rangle_{0,d}$ is computed via the following equality:
\ba
kt+\sum_{d=1}^{\infty}e^{dt}\langle{\cal O}_{h^{k-2-m}}{\cal O}_{h^{m-1}}\rangle_{0,d}=w({\cal O}_{h^{k-2-m}}{\cal O}_{h^{m-1}})_{0}(x(t)),
\ea
where $x=x(t)$ is the inverse of the mirror map $t=t(x)$.
We can easily see the following relation:
\ba
w({\cal O}_{h^{a}}{\cal O}_{h^{b}}|{\cal O}_{h})_{0,d}=d\cdot w({\cal O}_{h^{a}}{\cal O}_{h^{b}})_{0,d}.
\ea
Hence, the generating function $\tilde{L}^{k,k}_{m}(e^{x})$ is related to $w({\cal O}_{h^{k-2-m}}{\cal O}_{h^{m-1}})_{0}(x)$ in the following way:
\ba
\tilde{L}^{k,k}_{m}(e^{x}):=\frac{d}{dx}\frac{w({\cal O}_{h^{k-2-m}}{\cal O}_{h^{m-1}})_{0}(x)}{k}.
\ea
In particular, we have:
\ba
\tilde{L}^{k,k}_{1}(e^{x}):=\frac{dt}{dx}.
\ea
We also mention the following equality that was proved in \cite{Jin3}:
\ba
\tilde{L}^{k,k,}_{0}(e^{x}):=1+\sum_{d=1}^{\infty}\frac{w({\cal O}_{h^{k-2}}{\cal O}_{h^{-1}}|{\cal O}_{h})_{0,d}}{k}e^{dx}=\sum_{d=0}^{\infty}\frac{(kd)!}{(d!)^{k}}e^{dx}.
\ea
With this setup, the BCOV-Zinger formula\cite{BCOV, zinger} is given as follows.
\ba
F_{1}^{k,k,B}(x)&=&-\frac{1}{24}\left(\int_{M_{k}^{k}}c_{k-3}(T^{\prime}M_{k}^{k})\wedge h\right)x+\frac{1}{24}\chi(M_{k}^{k})\log\left(\tilde{L}_{0}^{k,k}(e^{x})\right)\no\\
                   &&-\frac{k-1}{48}\log(1-k^{k}e^{x})-\sum_{p=0}^{\frac{k-3}{2}}\frac{(k-1-2p)^2}{8}\log\left(\tilde{L}_{p}^{k.k}(e^{x})\right)\;\;(k: \mbox{odd}),\no\\
 F_{1}^{k,k,B}(x)&=&-\frac{1}{24}\left(\int_{M_{k}^{k}}c_{k-3}(T^{\prime}M_{k}^{k})\wedge h\right)x+\frac{1}{24}\chi(M_{k}^{k})\log\left(\tilde{L}_{0}^{k,k}(e^{x})\right)\no\\                  
&&-\frac{k-4}{48}\log(1-k^{k}e^{x})-\sum_{p=0}^{\frac{k-4}{2}}\frac{(k-2p)(k-2p-2)}{8}\log\left(\tilde{L}_{p}^{k.k}(e^{x})\right)\;\;(k:\mbox{even}),
\ea
where $\chi(M_{k}^{k})$ is Euler number of $M_{k}^{k}$ given by: 
\ba
\chi(M_{k}^{k})=\int_{M_{k}^{k}}c_{k-2}(T^{\prime}M_{k}^{k}).
\ea
In \cite{zinger}, Zinger proved the following theorem:
\begin{theorem}{\bf (Zinger \cite{zinger})}
\ba
F_{1}^{k,k,A}(t)&:=&-\frac{1}{24}\left(\int_{M_{k}^{k}}c_{k-3}(T^{\prime}M_{k}^{k})\wedge h\right)t+\sum_{d=1}^{\infty}\langle*\rangle_{1,d}e^{dt}\no\\
                  &=&F_{1}^{k,k,B}(x(t)),
\ea
where $\langle*\rangle_{1,d}$ is the Gromov-Witten invariant of $M_{k}^{k}$ of genus $1$ and of degree $d$ with no operator insertions, and $x=x(t)$ is inverse 
of the mirror map $t=t(x)$ given in (\ref{mirr}).
\label{zin}
\end{theorem}
At this stage, we revisit Conjecture \ref{main}.
In the $N=k$ case, $w(\prod_{p=2}^{N-2}({\cal O}_{h^{p}})^{m_{p}})_{1,d}$ is non-zero only if $m_{2}=m_{3}=\cdots=m_{N-2}=0$. Hence, $F_{1,vir.}^{k,k,B}$ has the 
following structure:
\ba
&&F_{1,vir.}^{k,k,B}(x^{1})\no\\
&=&-\frac{1}{24}\left( \int_{M_{k}^{k}}h\wedge c_{k-3}(T^{\prime}M_{k}^{k}) \right)x^{1}+\sum_{d=1}^{\infty}e^{dx^{1}}w(*)_{1,d},
\label{f1cyb}
\ea
where the symbol $*$ represents having no operator insertions. 
On the other hand, $F_{1}^{k,k,A}$ also has the same structure:
\ba
&&F_{1}^{k,k,A}(t^{1})\no\\
&=&-\frac{1}{24}\left( \int_{M_{k}^{k}}h\wedge c_{k-3}(T^{\prime}M_{k}^{k}) \right)t^{1}+\sum_{d=1}^{\infty}e^{dt^{1}}\langle*\rangle_{1,d}.
\ea
Therefore, we only need the mirror map that relates $x^{1}$ to $t^{1}$.
\ba
t^1(x^{1})=
x^{1}+\frac{1}{k}\sum_{d=1}^{\infty}e^{dx^{1}}w({\cal O}_{h^{k-3}}{\cal O}_{1}|*)_{0,d}.
\ea
Then Conjecture \ref{main} asserts the following equality:
\ba
F_{1}^{k,k,A}(t^{1})=F_{1.vir.}^{k,k,B}(x^{1}(t^{1})).
\ea
On the other hand, the mirror map which relates $x$ to $t$ in (\ref{mirr}) has the structure:
\ba
t(x)=x+\frac{1}{k}\sum_{d=1}^{\infty}e^{dx^{1}}w({\cal O}_{h^{k-3}}{\cal O}_{1})_{0,d}.
\ea  
By comparing  the r.h.s. of (\ref{multivir}) and (\ref{vir}), we can easily see that the following equality holds:
\ba
w({\cal O}_{h^{k-3}}{\cal O}_{1}|*)_{0,d}=w({\cal O}_{h^{k-3}}{\cal O}_{1})_{0,d}.
\ea 
In this way, we can identify $t^{1}$ and $x^{1}$ in Conjecture \ref{main} with $t$ and $x$ in Theorem \ref{zin}, respectively.
This naturally leads us to the following conjecture.
\begin{conj}
\ba
F_{1,vir.}^{k,k,B}(x)=F_{1}^{k,k,B}(x).
\ea
\label{jinzin}
\end{conj}
In the following, we present the results of numerical tests of this conjecture.

By Definition \ref{defw1}, $w(*)_{1,d}$ in (\ref{f1cyb}) is given by the sum of residues of integrands $f_{\Gamma}$ associated with $\Gamma\in \mbox{Graph}_{d}$:
\ba
w(*)_{1,d}=\sum_{\Gamma\in \mbox{\small Graph}_{d}}\mbox{Res}(f_{\Gamma}).
\ea 
We denote by $\mbox{Graph}^{(i)}_{d},\cdots, \mbox{Graph}^{(iv)}_{d}$ the sets of graphs of type (i),$\cdots$,(iv) of degree $d$, respectively.  We also denote by $\Gamma_{d}^{loop}$
and $\Gamma_{d}^{point}$ the unique graphs that belong to $\mbox{Graph}^{(ii)}_{d}$ and $\mbox{Graph}^{(iv)}_{d}$, respectively\footnote{Of course, we have to note that $\mbox{Graph}^{(ii)}_{1}=\emptyset$.}.   Then the above equality can be rewritten as follows.
\ba
w(*)_{1,d}&=&\sum_{\Gamma\in \mbox{\small Graph}^{(i)}_{d}}\mbox{Res}(f_{\Gamma})+\sum_{\Gamma\in \mbox{\small Graph}^{(ii)}_{d}}\mbox{Res}(f_{\Gamma})+
\sum_{\Gamma\in \mbox{\small Graph}^{(iii)}_{d}}\mbox{Res}(f_{\Gamma})+\sum_{\Gamma\in \mbox{\small Graph}^{(iv)}_{d}}\mbox{Res}(f_{\Gamma})\no\\
&=&\sum_{\Gamma\in \mbox{\small Graph}^{(i)}_{d}}\mbox{Res}(f_{\Gamma})+\mbox{Res}(f_{\Gamma_{d}^{loop}})+
\sum_{\Gamma\in \mbox{\small Graph}^{(iii)}_{d}}\mbox{Res}(f_{\Gamma})+\mbox{Res}(f_{\Gamma_{d}^{point}}),
\label{decom}
\ea 
At this stage, we introduce the following propositions.
\begin{prop}
For any positive $d$ and any $\Gamma\in \mbox{Graph}^{(iii)}_{d}$, $\mbox{Res}(f_{\Gamma})$ vanishes.
\end{prop}
{\it proof)}
Consider the integrand for a type (iii) graph and examine the residue integral in the $w$-variable:
\ban
&&\frac{1}{2\pi\sqrt{-1}}\oint_{C_{w}}dw\;\mbox{Sym}(\sigma)\frac{1}{24}\left(-\frac{k-1}{k}\frac{1}{w^k}-\frac{k+1}{k}\frac{1}{(z_{0})^{k}}\right)\frac{1}{(w-z_{0})^{2}}\frac{e_{k}(w,z_{0})}{kw}
F(z_{*,*})\\
&=&\frac{1}{2\pi\sqrt{-1}}\oint_{C_{w}}dw\;\mbox{Sym}(\sigma)\frac{1}{24}\left(-\frac{k-1}{k}\frac{1}{w^k}-\frac{k+1}{k}\frac{1}{(z_{0})^{k}}\right)\frac{1}{(w-z_{0})^{2}}\left(\prod_{j=1}^{k}((k-j)w+jz_{0})\right)
F(z_{*,*}).
\ean 
where $F(z_{*,*})$ is a rational function in $z_{0}$ and $z_{i,j}$'s, and $\frac{1}{2\pi\sqrt{-1}}\oint_{C_{w}}dw$ denotes the takeing the residue at $w=z_{0}$.
The result of the  above integration is given by:
\ban
&&\mbox{Sym}(\sigma)\frac{1}{24}\frac{d}{dw}\left.\left(\left(-\frac{k-1}{k}\frac{1}{w^k}-\frac{k+1}{k}\frac{1}{(z_{0})^{k}}\right)\left(\prod_{j=1}^{k}((k-j)w+jz_{0})\right)
F(z_{*,*}).\right)\right|_{w=z_{0}}
\ean
\ban
&=&\mbox{Sym}(\sigma)\frac{1}{24}\left((k-1)\frac{1}{(z_{0})^{k+1}}k^k(z_{0})^{k}-2\frac{1}{(z_{0})^{k}}\left(\sum_{j=1}^{k}(k-j)\right)k^{k-1}(z_{0})^{k-1}\right)F(z_{*,*})\\
&=&\mbox{Sym}(\sigma)\frac{1}{24}\frac{1}{z_{0}}\left((k-1)\cdot k^{k}-2\cdot\frac{k(k-1)}{2}\cdot k^{k-1}\right)F(z_{*,*})\\
&=&0.
\ean
Hence, the assertion of the proposition follows.  \hspace{2cm} $\Box$

\vspace{1cm}

\begin{prop}
The following equality holds.
\ba
\sum_{d=1}^{\infty}\left( \sum_{\Gamma\in \mbox{\small Graph}^{(i)}_{d}}\mbox{Res}(f_{\Gamma}) \right)e^{dx}=\frac{1}{24}\chi(M_{k}^{k})\log\left(\tilde{L}_{0}^{k,k}(e^{x})\right).
\label{prop2}
\ea
\end{prop}
{\it proof)}
Let us consider the integrand for a type (i) star graph obtained from the partition $\sigma=(d_{1},\cdots,d_{l(\sigma)})\in P_{d}$: 
\ba
f_{\Gamma}=&&\mbox{Sym}(\sigma)\frac{1}{24}\frac{c_{T}(z_{0})}{(kz_{0})^{l(\sigma)-1}}\left(\prod_{i=1}^{l(\sigma)}\frac{e_{k}(z_{0},z_{i,1})}{z_{i,1}-z_{0}}\right)
\left(\prod_{i=1}^{l(\sigma)}\left(\prod_{j=1}^{d_{i}-1}\frac{e_{k}(z_{i,j},z_{i,j+1})}{(2z_{i,j}-z_{i,j-1}-z_{i,j+1})kz_{i,j}}\right)\right)\no\\
&&\times \left(\frac{1}{(z_{0})^{k}}\right)\left(\prod_{i=1}^{l(\sigma)}\prod_{j=1}^{d_{i}}\frac{1}{(z_{i,j})^k}\right),\no\\
&&\hspace{9cm}(z_{i,0}=z_{0}\; (i=1,2,\cdots,l)).
\label{int2}
\ea
We first integrate over $z_{0}$. Collect factors involving $z_{0}$:
\ban
&&\frac{c_{T}(z_{0})}{(kz_{0})^{l(\sigma)-1}}\left(\prod_{i=1}^{l(\sigma)}\frac{e_{k}(z_{0},z_{i,1})}{z_{i,1}-z_{0}}\right)
\left(\prod_{i=1}^{l(\sigma)}\frac{1}{(2z_{i,1}-z_{0}-z_{i,2})}\right)\frac{1}{(z_{0})^k}\no\\
&=&kz_{0}\cdot c_{T}(z_{0})\left(\prod_{i=1}^{l(\sigma)}\frac{\prod_{j=1}^{k}\left((k-j)z_{0}+jz_{i,1}\right)}{z_{i,1}-z_{0}}\right)
\left(\prod_{i=1}^{l(\sigma)}\frac{1}{(2z_{i,1}-z_{0}-z_{i,2})}\right)\frac{1}{(z_{0})^k}\no\\
&=&\frac{\chi(M_{k}^{k})}{z_{0}}\left(\prod_{i=1}^{l(\sigma)}\frac{\prod_{j=1}^{k}\left((k-j)z_{0}+jz_{i,1}\right)}{z_{i,1}-z_{0}}\right)
\left(\prod_{i=1}^{l(\sigma)}\frac{1}{(2z_{i,1}-z_{0}-z_{i,2})}\right),
\ean  
where we used the equality $kz_{0}\cdot c_{T}(z_{0})=\chi(M_{k}^{k})(z_{0})^{k-1}$. Thus, after integrating $z_{0}$, the integrand becomes:
\ba
&&\mbox{Sym}(\sigma)\frac{\chi(M_{k}^{k})}{24}\left(\prod_{i=1}^{l(\sigma)}\frac{k!}{z_{i,1}(2z_{i,1}-z_{i,2})}\right)
\left(\prod_{i=1}^{l(\sigma)}\left(\prod_{j=2}^{d_{i}-1}\frac{e_{k}(z_{i,j},z_{i,j+1})}{(2z_{i,j}-z_{i,j-1}-z_{i,j+1})kz_{i,j}}\right)\right)\no\\
&&\times \left(\prod_{i=1}^{l(\sigma)}\prod_{j=2}^{d_{i}}\frac{1}{(z_{i,j})^k}\right).
\label{prep}
\ea
Next, we integrate the remaining variables. For this purpose, consider the residue integral:
\ban
&&C_{d}:=\frac{1}{(2\pi\sqrt{-1})^{d}}\oint_{C_{z_{1}}}dz_{1}\cdots\oint_{C_{z_{d}}}dz_{d}\;k!\cdot\frac{1}{z_{1}(2z_{1}-z_{2})}\left(\frac{\prod_{j=2}^{d}e_{k}(z_{j-1},z_{j})}{\prod_{i=2}^{d_{1}-1}kz_{i}(2z_{i}-z_{i-1}-z_{i+1})}\right)\left(\prod_{q=2}^{d}\frac{1}{(z_{q})^{k}}\right).
\ean 
Here, $\frac{1}{2\pi\sqrt{-1}}\oint_{C_{z_{1}}}dz_{1}$ denotes the operation of  taking residues at $z_{1}=0$ and $z_{1}=\frac{z_{2}}{2}$, 
$\frac{1}{2\pi\sqrt{-1}}\oint_{C_{z_{i}}}dz_{i}\;\;(i=2,\cdots, d-1)$
denotes the operation of taking residues at $z_{i}=0$ and $z_{i}=\frac{z_{i-1}+z_{i+1}}{2}$, and $\frac{1}{2\pi\sqrt{-1}}\oint_{C_{z_{d}}}dz_{d}$ denotes the operation of taking the residue at $z_{d}=0$. From (\ref{prep}), we can  see  that $\mbox{Res}(f_{\Gamma})$ is given by:
\ban
\mbox{Sym}(\sigma)\frac{\chi(M_{k}^{k})}{24}\prod_{i=1}^{l(\sigma)}C_{d_{i}}.
\ean
Then, we integrate over the $z_{1}$ variable of $C_{d}$. If we take the residue at $z_{1}=0$, we obtain the following contribution:
\ban
\frac{1}{(2\pi\sqrt{-1})^{d-1}}\oint_{C_{z_{2}}}dz_{2}\cdots\oint_{C_{z_{d}}}dz_{d}\;k!\cdot k!\cdot\frac{(-1)}{z_{2}(2z_{2}-z_{3})}\left(\frac{\prod_{j=3}^{d_{1}}e_{k}(z_{j-1},z_{j})}{\prod_{i=3}^{d-1}
kz_{i}(2z_{i}-z_{i-1}-z_{i+1})}\right)\left(\prod_{q=3}^{d}\frac{1}{(z_{q})^{k}}\right).
\ean
If we take the residue at $z_{1}=\frac{z_{2}}{2}$, we obtain:
\ban
\frac{1}{(2\pi\sqrt{-1})^{d-1}}\oint_{C_{z_{2}}}dz_{2}\cdots\oint_{C_{z_{d}}}dz_{d}\frac{(2k)!}{(2!)^{k}}\cdot \frac{1}{z_{2}(3z_{2}-2z_{3})}\left(\frac{\prod_{j=3}^{d_{1}}e_{k}(z_{j-1},z_{j})}{
\prod_{i=3}^{d-1}kz_{i}(2z_{i}-z_{i-1}-z_{i+1})}\right)\left(\prod_{q=3}^{d}\frac{1}{(z_{q})^{k}}\right).
\ean
We then introduce a set of ordered partitions of the positive integer $d$.
\ban
OP_{d}:=\{\;(\tau=(d_{1},d_{2},\cdots,d_{l(\tau)})\;|\;d_{j}\geq 1\;(j=1,2,\cdots,l(\tau)),\;\sum_{j=1}^{l(\tau)}d_{j}=d\;\}.
\ean
Continuing with the computations, we can express $C_{d}$ explicitly as:
\ban
C_{d}&=&\sum_{\tau=(d_{1},\cdots,d_{l(\tau)})\in OP_{d}}(-1)^{l(\tau)-1}\prod_{j=1}^{l(\tau)}\frac{(kd_{j})!}{(d_{j}!)^{k}}\\
      &=& \sum_{\sigma=(d_{1},\cdots,d_{l(\sigma)})\in P_{d}}(-1)^{l(\sigma)-1}\frac{l(\sigma)!}{\prod_{j=1}^{d}\mbox{mul}(j;\sigma)!}\prod_{j=1}^{l(\sigma)}\frac{(kd_{j})!}{(d_{j}!)^{k}},
\ean
where $\displaystyle{\frac{l(\sigma)!}{\prod_{j=1}^{d}\mbox{mul}(j;\sigma)!}}$ represents the number of ordered partitions that have the same elements ($d_{j}$) as those of $\sigma\in P_{d}$.
By using the well-known combinatorial identity:
\ban
1+\sum_{d=1}^{\infty}\left( \sum_{\sigma\in P_{d}}(-1)^{l(\sigma)}\frac{l(\sigma)!}{\prod_{j=1}^{d}\mbox{mul}(\sigma;j)!}\left(\prod_{j=1}^{l(\sigma)}b_{d_{j}}\right)\right)e^{dx}
=\left(1+\sum_{d=1}^{\infty}b_{d}e^{dx}\right)^{-1},
\ean
we obtain the equality:
\ba
1-\sum_{d=1}^{\infty}C_{d}e^{dx}=\left(1+\sum_{d=1}^{\infty}\frac{(kd)!}{(d!)^{k}}e^{dx}\right)^{-1}\left(=\frac{1}{\tilde{L}_{0}^{k,k}(e^{x})}\right).
\label{key1}
\ea
With these preparations, the l.h.s. of (\ref{prop2}) is evaluated as follows.
\ba
&&\sum_{d=1}^{\infty}\left( \sum_{\Gamma\in \mbox{\small Graph}^{(i)}_{d}}\mbox{Res}(f_{\Gamma}) \right)e^{dx}\no\\
&=&\sum_{d=1}^{\infty}\left( \sum_{\sigma\in P_{d}}\frac{(l(\sigma)-1)!}{\prod_{j=1}^{d}\mbox{mul}(\sigma;j)!}\frac{\chi(M_{k}^{k})}{24}\left(\prod_{j=1}^{l(\sigma)}C_{d_{j}}\right)\right)e^{dx}\no\\
&=&\frac{\chi(M_{k}^{k})}{24}\sum_{d=1}^{\infty}\left( \sum_{\sigma\in P_{d}}\frac{(l(\sigma)-1)!}{\prod_{j=1}^{d}\mbox{mul}(\sigma;j)!}\left(\prod_{j=1}^{l(\sigma)}C_{d_{j}}\right)\right)e^{dx}\no\\
&=&-\frac{\chi(M_{k}^{k})}{24}\log\left(1-\sum_{d=1}^{\infty}C_{d}e^{dx}\right).
\label{pf}
\ea
In going from the third line to the fourth line, we used another well-known combinatorial identity:
\ban
\sum_{d=1}^{\infty}\left( \sum_{\sigma\in P_{d}}\frac{(l(\sigma)-1)!}{\prod_{j=1}^{d}\mbox{mul}(\sigma;j)!}\left(\prod_{j=1}^{l(\sigma)}b_{d_{j}}\right)\right)e^{dx}
=-\log\left(1-\sum_{d=1}^{\infty}b_{d}e^{dx}\right).
\ean
On the other hand, by taking the logarithm of both sides of (\ref{key1}), we obtain,
\ba
\log(1-\sum_{d=1}^{\infty}C_{d}e^{dx})=-\log\left(\tilde{L}_{0}^{k,k}(e^{x})\right).
\label{key2}
\ea
By combining (\ref{pf}) with (\ref{key2}), we obtain the assertion of the proposition. \hspace{2cm}$\Box$

\vspace{1cm}

Let us return to the discussion of (\ref{decom}). $\mbox{Res}(f_{\Gamma_{d}^{point}})$ is immediately computed from (\ref{point}).
\ba
\mbox{Res}(f_{\Gamma_{d}^{point}})&=&\frac{1}{24}R_{k,k}(d)k^{kd}\no\\
&=&\frac{1}{24}\left(\frac{(k-1)}{2}\frac{1}{d}-\frac{k^2-1}{k}\frac{1}{d^2}\right)k^{kd}.
\label{point2}
\ea
Hence we obtain the following propostion\footnote{The expression in Proposition \ref{yukawa} is closely related to the logarithm of the B-model Yukawa coupling, or the singular locus of the mirror manifold of the quintic 3-fold. However, the reason for this connection is not clearly understood within our formalism.}.
\begin{prop}
\ba
\sum_{d=1}^{\infty}\mbox{Res}(f_{\Gamma_{d}^{point}})e^{dx}=-\frac{k-1}{48}\log(1-k^ke^{x})+\frac{k^2-1}{24k}\int_{-\infty}^{x}\log(1-k^ke^{s})ds.
\ea
\label{yukawa}
\end{prop}
By combining these results, we obtain,
\ba
F_{1,virt.}^{k,k,B}(x)&=&-\frac{1}{24}\left(\int_{M_{k}^{k}}c_{k-3}(T^{\prime}M_{k}^{k})\wedge h\right)x+\frac{1}{24}\chi(M_{k}^{k})\log\left(L_{0}^{k,k}(e^{x})\right)\no\\
                   &&+\sum_{d=2}^{\infty} \mbox{Res}(f_{\Gamma_{d}^{loop}})e^{dx} -\frac{k-1}{48}\log(1-k^ke^{x})+\frac{k^2-1}{24k}\int_{-\infty}^{x}\log(1-k^ke^{s})ds.
\label{tbvir}
\ea
{\bf This is another representation of the BCOV-Zinger formula obtained from Conjecture \ref{jinzin}.}
Therefore, Conjecture \ref{jinzin} is restated as the following identities.
\begin{conj}
\ba
&&\sum_{d=2}^{\infty} \mbox{Res}(f_{\Gamma_{d}^{loop}})e^{dx}+\frac{k^2-1}{24k}\int_{-\infty}^{x}\log(1-k^ke^{s})ds\no\\
&=&-\sum_{p=0}^{\frac{k-3}{2}}\frac{(k-1-2p)^2}{8}\log\left(\tilde{L}_{p}^{k.k}(e^{x})\right)\;\;(k:\mbox{odd}),\no\\
&&\sum_{d=2}^{\infty} \mbox{Res}(f_{\Gamma_{d}^{loop}})e^{dx}-\frac{1}{16}\log(1-k^ke^{x})+\frac{k^2-1}{24k}\int_{-\infty}^{x}\log(1-k^ke^{s})ds\no\\
&=&-\sum_{p=0}^{\frac{k-4}{2}}\frac{(k-2p)(k-2p-2)}{8}\log\left(\tilde{L}_{p}^{k.k}(e^{x})\right)\;\;(k:\mbox{even}).
\label{idconj}
\ea
\label{loopconj}
\end{conj}
We numerically confirmed these identities for $k=4,5,6,7,8$ up to $d=5$ \footnote{Fundamentally, it is possible to prove Conjecture \ref{loopconj} by direct computation of the residue integral. However, due to its complexity, we leave the proof for future work.}. Below, we present here numerical data up to $d=5$ in the $k=5$ (quintic 3-fold) case.

The generating function of loop amplitudes is given by:
\ba
&&\sum_{d=2}^{\infty} \mbox{Res}(f_{\Gamma_{d}^{loop}})e^{dx}\no\\
&&=-\frac{1174875}{4} e^{2x}-\frac{6913090625}{9} e^{3x}-\frac{31054165371875}{16} e^{4x}-5008379074144375 e^{5x}-\cdots.
\ea 
We also have:
\ba
&&\frac{1}{5}\int_{-\infty}^{x}\log(1-5^5e^{s})ds\no\\
&&=-625 e^x-\frac{1953125}{4} e^{2x}-\frac{6103515625}{9} x^3-\frac{19073486328125}{16} e^{4x}-2384185791015625 e^{5x}-\cdots.
\ea 
Then, the l.h.s. of (\ref{idconj}) turns out to be: 
\ba
&&\sum_{d=2}^{\infty} \mbox{Res}(f_{\Gamma_{d}^{loop}})e^{dx}+\frac{1}{5}\int_{-\infty}^{x}\log(1-5^5e^{s})ds\no\\
&&=-625e^x-782000 e^{2x}-\frac{4338868750}{3} e^{3x}-3132978231250 e^{4x}-7392564865160000 e^{5x}-\cdots.
\ea
 On the other hand, $\log\left(\tilde{L}_{0}^{5.5}(e^{x})\right)$ and $\log\left(\tilde{L}_{1}^{5.5}(e^{x})\right)$ are given as follows.
 \ba
 &&\log\left(\tilde{L}_{0}^{5.5}(e^{x})\right)\no\\
 &&=120 e^{x}+106200 e^{2x}+155136000 e^{3x}+280511415000 e^{4x}+571399451565120 e^{5x}+\cdots,\no\\
  &&\log\left(\tilde{L}_{1}^{5.5}(e^{x})\right)\no\\
 &&=770e^{x} + 1139200e^{2x} + \frac{6816105500}{3}e^{3x} + 5143910802500e^{4x} + 12499531924059520e^{5x}+\cdots.
 \ea
 And the r.h.s. of (\ref{idconj}) becomes:
 \ba
 &&-2\log\left(\tilde{L}_{0}^{5.5}(e^{x})\right)-\frac{1}{2}\log\left(\tilde{L}_{1}^{5.5}(e^{x})\right)\no\\
 &&=-625e^x-782000 e^{2x}-\frac{4338868750}{3} e^{3x}-3132978231250 e^{4x}-7392564865160000 e^{5x}-\cdots.
 \ea
 In this way, the identity (\ref{idconj}) in the $k=5$ case is numerically confirmed up to $d=5$.

\newpage

\appendix

\section{Tables for N=5 cases}

\begin{table}[H]
\centering
\caption{$M_5^1$}
\begin{tabular}{|c|l|l|l|l|l|}
\hline
d&(a,b)&$N^0_{d,a,b}$ & $N^1_{d,a,b}$ &$\frac{2d-1}{12}N^0_{d,a,b}+N^1_{d,a,b}$&$w_{a,b}$  \\
\hline
1&(0,2) & 1 &$-\frac{1}{12}$&0&$-\frac{7}{12}$\\ \hline
1&(2,1) & 1 &$-\frac{1}{12}$&0&$-\frac{5}{6}$\\ \hline
1&(4,0) & 2 &$-\frac{1}{6}$&0&$-\frac{7}{6}$\\ \hline
2&(0,4)& 0 & $0$& 0 &$-\frac{76}{3}$\\ \hline
2&(2,3)& 1 & $-\frac{1}{4}$& 0&$-\frac{853}{12}$ \\ \hline
2& (4,2)& 4 &$ -1$            &0&$-198$ \\ \hline
2& (6,1)& 18&$-\frac{9}{2}$& 0&$-\frac{1097}{2}$\\ \hline
2& (8,0) & 92&$-23$           & 0 &$-\frac{4541}{3}$ \\ \hline
3&(0,6)& 1 & $-\frac{5}{12}$& 0 &$-\frac{19959}{4}$\\ \hline
3&(2,5)& 5 & $-\frac{25}{12}$& 0&$-\frac{62338}{3}$ \\ \hline
3& (4,4)& 30&$ -\frac{25}{2}$&0&$-\frac{516827}{6}$ \\ \hline
3& (6,3)& 190&$-\frac{469}{6}$& 1&$-\frac{1068442}{3}$\\ \hline
3& (8,2) & 1312&$-\frac{1598}{3} $& 14&$-\frac{4408330}{3}$ \\ \hline
3& (10,1) & 9864&$-3960$ & 150 &$-\frac{18159922}{3}$\\ \hline
3& (12,0) & 80160&$-31900$ & 1500 &$-\frac{74719852}{3}$\\ \hline
4&(0,8)& 4 & $-\frac{4}{3}$& 1&$-\frac{7111330}{3}$ \\ \hline
4&(2,7)& 58 & $-\frac{179}{6}$& 4&$-\frac{26141813}{2}$ \\ \hline
4& (4,6)& 480&$ -248 $&32 &$-71830274$\\ \hline
4& (6,5)& 4000&$-\frac{6070}{3}$& 310&$-\frac{1182256279}{3}$\\ \hline
4& (8,4) & 35104&$-\frac{51772}{3} $& 3220&$-2159333004$ \\ \hline
4& (10,3) &327888&$-156594$  & 34674&$-\frac{35458691818}{3}$ \\ \hline
4& (12,2) & 3259680&$-1515824$  & 385656&$-\frac{193936379144}{3}$ \\ \hline
4& (14,1) & 34382544&$-15620216$& 4436268 &$-353359995764$\\ \hline
4& (16,0) &383306880&$-170763640$ & 52832040&$-1930689790136$ \\ \hline
5&(0,10)& 105 & $-\frac{147}{4}$& 42 &$-\frac{8363354113}{4}$ \\ \hline
5&(2,9)& 1265 & $-\frac{2379}{4} $& 354 &$-\frac{28682135389}{2}$\\ \hline
5& (4,8)& 13354&$ -\frac{13047}{2} $&3492 &$-\frac{196198477325}{2}$\\ \hline
5& (6,7)&139098&$-\frac{132549}{2}$& 38049 &$-\frac{2010681907978}{3}$\\ \hline
5& (8,6) & 1492616&$-677808 $& 441654 &$-\frac{13724961403006}{3}$ \\ \hline
5& (10,5) &16744080&$-7179606$   & 5378454 &$-\frac{93619004917238}{3}$\\ \hline
5& (12,4) & 197240400&$-79637976$ & 68292324 &$-212735629674372$ \\ \hline
5& (14,3) & 2440235712&$-928521900$& 901654884 &$-\frac{4348697671027760}{3}$\\ \hline
5& (16,2) &31658432256&$-11385660384$ & 12358163808 &$-9873859605646752$\\ \hline
5& (18,1) &429750191232&$-146713008096$ & 175599635328 &$-\frac{201722432909390752}{3}$\\ \hline
5& (20,0) &6089786376960&$-1984020394752$  & 2583319387968 &$-\frac{1373530281059327936}{3}$\\ \hline
\end{tabular}
\end{table}

\begin{table}[H]
\centering
\caption{$M_5^2$}
\begin{tabular}{|c|l|l|l|l|l|}
\hline
d&(a,b)&$N^0_{d,a,b}$ & $N^1_{d,a,b}$ &$\frac{3d-2}{24}N^0_{d,a,b}+N^1_{d,a,b}$&$w_{a,b}$  \\
\hline
1&(1,1) & 4 &$-\frac{1}{6}$&0&$-\frac{13}{6}$\\ \hline
1&(3,0) & 8 &$-\frac{1}{3}$&0&$-3$\\ \hline
2&(0,3)& 8 & $-\frac{4}{3}$& 0&$-\frac{287}{3}$ \\ \hline
2& (2,2)&16 &$ -\frac{8}{3}$ &0&$-264$ \\ \hline
2& (4,1)& 64&$-\frac{32}{3}$& 0&$-\frac{2174}{3}$\\ \hline
2& (6,0) & 320&$-\frac{160}{3}$ & 0 &$-\frac{5956}{3}$ \\ \hline
3&(1,4)& 64 & $-\frac{56}{3}$& 0 &$-\frac{104500}{3}$\\ \hline
3&(3,3)& 320 & $-\frac{280}{3}$& 0&$-\frac{429196}{3}$ \\ \hline
3& (5,2)&2048&$ -\frac{1792}{3}$&0&$-\frac{1759552}{3}$ \\ \hline
3& (7,1)& 15104&$-\frac{13216}{3}$& 0&$-\frac{7209584}{3}$\\ \hline
3& (9,0) & 123904&$-\frac{108416}{3} $& 0&$-\frac{29527616}{3}$ \\ \hline
4& (0,6)& 384&$ -160 $&0 &$-\frac{18667312}{3}$\\ \hline
4& (2,5)& 2560&$-\frac{3200}{3}$& 0&$-\frac{101879272}{3}$\\ \hline
4& (4,4) & 18944&$-\frac{22912}{3}$& 256&$-\frac{555449168}{3}$ \\ \hline
4& (6,3) &163840&$-\frac{194048}{3}$ & 3584&$-\frac{3026251616}{3}$ \\ \hline
4& (8,2) &1583104&$-\frac{1849856}{3}$  &43008&$-\frac{16485590720}{3}$ \\ \hline
4& (10,1) & 16687104&$-6440960$&512000 &$-\frac{89806527616}{3}$\\ \hline
4& (12,0) &189358080&$-72652800$ & 6246400&$-163085218816$ \\ \hline
5& (1,7)&27136&$-\frac{41792 }{3}$& 768 &$-\frac{28726121392}{3}$\\ \hline
5& (3,6) & 229376&$-\frac{331264}{3}$&13824 &$-\frac{195282001984}{3}$ \\ \hline
5& (5,5) &2232320&$-\frac{3049984}{3}  $ &192512 &$-\frac{1326874482304}{3}$\\ \hline
5& (7,4) &24391680&$-10660352$ & 2551808 &$-\frac{9013280450048}{3}$ \\ \hline
5& (9,3) &291545088&$-123583488$ & 34336768 &$-\frac{61226330115584}{3}$\\ \hline
5& (11,2) &3750199296&$-1553444864$  &477913088 &$-138652119786496$\\ \hline
5& (13,1) &51384877056&$-20917362688$  & 6916112384 &$-\frac{2826429058966016}{3}$\\ \hline
5& (15,0) &744875950080&$-299359264768$  & 104115208192 &$-\frac{19209989184830464}{3}$\\ \hline
\end{tabular}
\end{table}

\begin{table}[H]
\centering
\caption{$M_5^3$}
\begin{tabular}{|c|l|l|l|l|l|}
\hline
d&(a,b)&$N^0_{d,a,b}$ & $N^1_{d,a,b}$ &$\frac{d-1}{12}N^0_{d,a,b}+N^1_{d,a,b}$&$w_{a,b}$  \\
\hline
1&(0,1) & $18$ &0&0&$-\frac{21}{2}$\\ \hline
1&(2,0) & $45$ &0&0&$-\frac{27}{2}$\\ \hline
2&(0,2)&$ 54$  & $-\frac{9}{2}$& 0&$-\frac{2187}{2}$ \\ \hline
2& (2,1)&$378$ &$ -\frac{63}{2}$ &0&$-2862$ \\ \hline
2& (4,0)&$2187$&$-\frac{729}{4}$& 0&$-\frac{30501}{4}$\\ \hline
3&(0,3)& $648$ &$ -81$& $27$&$-299943$ \\ \hline
3& (2,2)&$7452$&$-1161$&$81$&$-1188027$ \\ \hline
3& (4,1)&$65610$ &$-10449$& $486$&$-\frac{9537669}{2}$\\ \hline
3& (6,0) & $623295$&$-\frac{200475}{2} $&$3645$&$-19201644$ \\ \hline
4& (0,4) & $15552$&$-1701$&$2187$&$-\frac{279086715}{2}$ \\ \hline
4& (2,3) &$248832$&$-48357$ & $13851$&$-740281275$ \\ \hline
4& (4,2) &$2991816$&$-616734$  &$131220$&$-3968742582$ \\ \hline
4& (6,1) & $37161504$&$-7846956$&$1443420$ &$-21341675475$\\ \hline
4& (8,0) &$491956902$&$-\frac{211336371}{2}$ &$17321040$&$-\frac{229871126583}{2}$ \\ \hline
5& (0,5) &$583200$&$972$ &$195372$ &$-92893454856$\\ \hline
5& (2,4) &$11955600$&$-2005236$ &1979964 &$-617546315223$ \\ \hline
5& (4,3) &$183760488$&$-35803377$ &$25450119$ &$-\frac{8281495651131}{2}$\\ \hline
5& (6,2) &$2838367332$&$-584014293$  &$362108151$ &$-27840191130297$\\ \hline
5& (8,1) &$45746559378$&$-9717064074$  &$5531789052$ &$-187479083534526$\\ \hline
5& (10,0) &$776682421065$&$-169540839261$  & $89353301094$ &$-1263878784214992$\\ \hline
\end{tabular}
\end{table}

\begin{table}[H]
\centering
\caption{$M_5^4$}
\begin{tabular}{|c|l|l|l|l|l|}
\hline
d&(a,b)&$N^0_{d,a,b}$ & $N^1_{d,a,b}$ &$\frac{d-2}{24}N^0_{d,a,b}+N^1_{d,a,b}$&$w_{a,b}$  \\
\hline
1&(1,0) & $320$ &$\frac{40}{3}$&0&$-\frac{344}{3}$\\ \hline
2&(0,1)&$3888$  & $0$& 0&$-\frac{84848}{3}$ \\ \hline
2& (2,0)&$27200$ &$0$ &0&$-\frac{222080}{3}$ \\ \hline
3&(1,1)& $672768$ &$-24192$& $3840$&$-32895232$ \\ \hline
3& (3,0)&$8388608$&$-\frac{971776}{3}$&$25600$&$-\frac{390036992}{3}$ \\ \hline
4& (0,2) &$18323712$&$-861696$  &$665280$&$-13842672128$ \\ \hline
4& (2,1) & $284802048$&$-14229504$&$9504000$ &$-\frac{221036603392}{3}$\\ \hline
4& (4,0) &$5100273664$&$-\frac{905199616}{3}$ &$123289600$&$-\frac{1171415220224}{3}$ \\ \hline
5& (1,2) &$9830744064$&$-148801536$  &$1080041472$ &$-\frac{113572126965760}{3}$\\ \hline
5& (3,1) &$206561083392$&$-4948770816 $  &$20871364608$ &$-252561662754816$\\ \hline
5& (5,0) &$4821100789760$&$-202932748288$  & $399704850432$ &$-\frac{5041038692581376}{3}$\\ \hline
\end{tabular}
\end{table}

\vspace{10cm}

\newpage
\noindent
{\bf Conflict of interest statement}\\ 
The authors have no competing interests to declare that are relevant to the content of this article.

\end{document}